\documentclass[hidelinks,11pt]{article}

\newcommand{\bol}{\boldsymbol}

\newcommand{\ner}{\bol{r}}

\newcommand{\de}{\,\mathrm{d}}                               
\newcommand{\e}{\operatorname{e}}                               
                     
\newcommand{\inc}{\mathrm{inc}}

\newcommand{\andtext}{\quad\mbox{and}\quad}

\newcommand{\p}{\partial}

\newcommand{\real}{\mathrm{Re}\,}    
                                   
\newcommand{\imag}{\mathrm{Im}\,}

\newcommand{\lf}{\left}
\newcommand{\rg}{\right}

\newcommand{\Z}{\mathbb{Z}}     
\newcommand{\R}{\mathbb{R}}       
\newcommand{\C}{\mathbb{C}}       
\newcommand{\N}{\mathbb{N}}

\newcommand{\nor}{{\bf n}}

\usepackage{amsmath}
\usepackage{amsmath}
\usepackage{multirow}
\usepackage{amsmath}
\usepackage{amsfonts}
\usepackage{amsmath}
\usepackage{amssymb}  
\usepackage{graphicx}
\usepackage{caption}
\usepackage{mathrsfs}
\usepackage{upgreek}
\usepackage{amsthm}
\usepackage{subfig}
\usepackage{booktabs}
\usepackage{authblk}
\usepackage{xcolor}
\usepackage{cite}

\usepackage[]{algorithm2e}

\newtheorem{theorem}{Theorem}[section]

\newtheorem{proposition}[theorem]{Proposition}

\newtheorem{remark}[theorem]{Remark}

 \usepackage{xcolor}
 \usepackage[normalem]{ulem}

\title{Windowed Green function method for wave scattering by periodic arrays of 2D obstacles\thanks{This work was supported by FONDECYT (Fondo Nacional de Desarrollo Científico y Tecnológico), Chile, Grant Number 11181032, and INRIA Chile-Pontificia Universidad Católica “The Bridge" exchange program.}}

 \author[1]{Thomas Strauszer-Caussade}
\affil[1]{\small{Institute for Mathematical and Computational Engineering, Pontificia Universidad Cat\'olica de Chile, Santiago, Chile}}
\author[2]{Luiz M. Faria}
 \affil[2]{\small{POEMS, CNRS, INRIA, ENSTA Paris, Institut Polytechnique de Paris, Palaiseau, France}}
\author[3]{Agust\'in Fernandez-Lado}
 \affil[3]{\small{Intel Corporation, Aloha, Oregon, USA}}
 \author[4]{\\ Carlos P\'erez-Arancibia\thanks{Corresponding author: c.a.perezarancibia@utwente.nl}}
 \affil[4]{\small{Department of Applied Mathematics and MESA+ Institute, University of Twente, Enschede,  The Netherlands}}

\date{\today}

\begin{document}
\maketitle


\begin{abstract}

This paper introduces a novel boundary integral equation (BIE) method for the numerical solution of problems of planewave scattering by periodic line arrays of two-dimensional penetrable obstacles. Our approach is built upon a direct BIE formulation that leverages the simplicity of the free-space Green function but in turn entails evaluation of integrals over the unit-cell boundaries. Such integrals are here treated via the  window Green function method. The windowing approximation together with a finite-rank operator correction---used to properly impose the Rayleigh radiation condition---yield a robust second-kind BIE that produces super-algebraically convergent solutions throughout the spectrum, including at the challenging Rayleigh-Wood anomalies. The corrected windowed BIE can be discretized by means of off-the-shelf Nystr\"om and boundary element methods, and it leads to linear systems suitable for iterative linear-algebra solvers as well as standard fast matrix-vector product algorithms. A variety of numerical examples demonstrate the accuracy and robustness of the proposed methodology.\\


\end{abstract}

\maketitle 

\section{Introduction}

This paper presents a novel windowed Green function boundary integral equation (BIE) method for the numerical solution of problems of time-harmonic electromagnetic planewave scattering by infinite periodic arrays of penetrable obstacles in two spatial dimensions (although the proposed methodology can also be applied to acoustics). Problems of this type often arise in a number of application areas that greatly benefit from accurate and efficient numerical computations such as, for instance, photonic crystal modeling~\cite{joannopoulos2008} and inverse design of metasurfaces~\cite{yu2014flat,xie2014wavefront} whereby the so-called \emph{locally periodic approximation} is used to deal with scattering by large aperiodic structures by decomposing it in a finite number of unit-cell periodic problems~\cite{pestourie2018inverse,perez2018sideways,lin2019topology}.

Classical BIE formulations for scattering by periodic media rely on the  quasi-periodic Green function~\cite{linton1998green}. As is well known, standard spatial and spectral representations of the quasi-periodic Green function give rise to infinite series that (a) converge slowly depending on the relative location of the source and target points and, in addition, (b) cease to exist at the so-called Rayleigh-Wood (RW) anomalies (i.e., when at least one scattered/transmitted mode propagates in the direction parallel to the array axis). Several analytical techniques have been proposed to tackle the former problem including most notably Ewald's method~(see~\cite{linton2010lattice,linton1998green} for a thorough review on the subject). A strikingly simple method that also addresses the aforementioned slow convergence issue is developed in~\cite{MonroJr:2008te,bruno2016superalgebraically}, which relies on a smooth windowed sum approximation of the spatial series representation of the Green function.  Away from RW anomalies, this approach  achieves super-algebraically fast convergence as the truncation radius increases. In view of the fact that the quasi-periodic Green function itself does not exists at RW anomalies, all the aforementioned approaches simply break down at these singular configurations (although, as shown in~\cite[Fig. 1.3]{fernandez2020wave}, Ewald's method produces accurate solutions at almost machine precision ``distance" from RW anomalies). 

Improving on the windowed summation approach,~\cite{brunodelourme} and subsequent related work~\cite{bruno2017rapidly,bruno2017three} introduce the quasi-periodic shifted Green function. BIE solvers that leverage this modified Green function~\cite{bruno2019shifted,nicholls2020sweeping,arancibia2019domain} exhibit super-algebraic convergence away from RW anomalies and algebraic but arbitrarily high-order convergence at and around RW anomalies, at the cost of $n$-tupling the number of function evaluations where $n$ is the numbers of ``shifts" utilized in the approximation. Recent developments in this direction~\cite{fernandez2020wave,bruno2020evaluation} present a general methodology based on hybrid spatial/spectral Green function representations and the Woodbury-Sherman-Morrison formula that makes classical approaches such as Laplace-type integral and Ewald's methods, as well as the shifted Green function approach, applicable and robust at and around RW-anomaly configurations.

Yet another class of BIE methods aims at bypassing the use of the problematic quasi-periodic Green function. The approach adopted  in~\cite{barnett2011new} (see also~\cite{gillman2013fast})  recasts the quasi-periodic problem as a formally second-kind indirect BIE formulation involving free-space Green function kernels and integrals along the infinite boundaries of the (unbounded) unit cell domain. Leveraging the exponential decay of the boundary integrands in spectral form, the resulting BIE system is effectively reduced to a bounded hybrid spatial-spectral computational domain where standard Nystr\"om discretizations can be applied for its numerical solution. Although this method does not make use of the quasi-periodic Green function, it involves evaluation of cumbersome Sommerfeld-type integrals that need to be painstakingly modified in the presence of RW anomalies (when a pole at origin on the Sommerfeld integration contour needs to be accounted for by suitably deforming the contour and adding the corresponding residue contribution). Building up on this work, a periodizing scheme akin to the method of fundamental solutions is developed in~\cite{cho2015robust} and subsequent contributions~\cite{lai2015fast}. This approach only entails evaluations of free-space Green function kernels in spatial form and it appears immune to the presence of RW anomalies. However, the so-called proxy (equivalent) sources employed by this scheme to enforce the quasi-periodicity condition, give rise to relatively small but ill-conditioned subsystems that are treated by Schur complements and direct linear algebra solvers, hence hindering the straightforward applicability of GMRES and fast algorithms to perform matrix-vector product operations (such as the fast multipole method~\cite{rokhlin1990rapid}, for instance).

Here, we present a method that falls under the latter class of BIE methods. Our approach amounts to an extension of the windowed Green function (WGF) method for layer media scattering a waveguide problems~\cite{bruno2016windowed,perez2017windowed,bruno2017windowed,bruno2017wg}, to quasi-periodic scattering problems. Inspired by~\cite{bruno2016windowed} we pursue a direct BIE formulation derived from a Green's representation formula of the scattered field within the unbounded unit-cell domain, which uses the free-space Green function instead of  the problem-specific (quasi-periodic) Green function. The quasi-periodicity condition is then readily incorporated into our formulation by exploiting the direct linear relationship between the scattered-field traces on the left- and right-hand side unit-cell walls. The transmission conditions on the penetrable boundaries of the obstacles are imposed through Kress-Roach/M\"uller's approach~\cite{kress1978transmission,muller2013foundations} which leads to weakly-singular integral operators. As in~\cite{barnett2011new}, we hence obtain a formally second-kind BIE system given in terms of free-space Green function kernels and boundary integrals over the unbounded unit-cell boundaries. Indeed, prior to truncation, both formulations entail evaluation of the same weakly-singular integral operators. The main difference between the two approaches lies in the truncation strategy. Instead of resorting to spatial-spectral representations of the integral operators, we work with integral operators in pure spatial form hence enabling the use of off-the-shelf BIE methods and fast algorithms. We do so by truncating the oscillatory integrals over the unbounded unit-cell walls using a smooth window function that multiplies the free-space Green function kernels.  When applied to the traces of the (radiative) scattered field, the resulting windowed BIE operators spawn small errors that decay super-algebraically fast as the support of the window function increases.  As it turns out, however, at certain frequency ranges which include RW-anomaly configurations, the naive windowing approximation of the BIE operators leads to a BIE system that fails to account for the radiation condition. In order to properly enforce it, we then propose a corrected windowed BIE that produces accurate solutions throughout the entire spectrum, including at and around the challenging RW-anomaly configurations. The corrected windowed BIE is Fredholm of the second-kind and upon discretization it leads to systems of equations that can be efficiently solved by iterative linear algebra solvers (i.e., GMRES) which can be further accelerated by means of fast methods. The proposed methodology exhibits super-algebraic convergence as the window size increases.

The paper is organized as follows. Section~\ref{sec:Inte_Eq_For} describes the problem under consideration and summarizes some important facts of the problem that will be utilized in the following sections. Section~\ref{sec:green_representation} presents the derivation of the non-standard Green's representation formula on which our direct BIE formulation is based on. Section~\ref{sec:param_op} introduces the notation as well as the main properties of the layer potentials and BIE operators.  The direct BIE formulation of the problem is derived in Section~\ref{sec:BIE} while the naive windowed BIE is motivated and presented in Section~\ref{sec:WGF}. A series of numerical experiments designed to examine the accuracy of the naive windowed BIE at and around a RW-anomaly configurations is shown in Section~\ref{sec:illus_num_exam}. The corrected windowed BIE formulation is then developed in Section~\ref{sec:corrected_BIE}. A variety of the numerical examples are presented in Section~\ref{sec:num_exam}. Finally, Section~\ref{sec:conclusions}, presents the conclusions and future work.

\section{Preliminaries}\label{sec:Inte_Eq_For}
This paper deals with problems of time-harmonic electromagnetic scattering by infinite periodic arrays of penetrable obstacles in two dimensions for which we adopt the time convention $\e^{-i\omega t}$ where $t>0$ is time and $\omega>0$ is the angular frequency. In detail, letting $\theta^\inc\in \lf[-\frac{\pi}2,\frac{\pi}{2}\rg]$ denote the angle of incidence measured with respect to the negative $y$-axis, we consider the scattering and transmission of a planewave 
\begin{equation}\label{eq:planewave}
u^{\rm inc}(x,y)= \e^{i\alpha x-i\beta y},\quad \alpha = k_1\sin\theta^\inc,\quad\beta = k_1\cos\theta^\inc
\end{equation}
by an $L$-periodic  array of the form $D =\bigcup_{n\in\Z} \Omega_n$ with $\Omega_n = \{(x,y)\in\R^2:(x-nL,y)\in\Omega\}$ where $\Omega\subset\R^2$ is an open and bounded domain. For the sake of ease of exposition the boundary $\Gamma_1:=\p\Omega$ of $\Omega = \Omega_0$ is assumed of class $C^2$ and given by
$$\Gamma_1:=\{\ner\in\R^2:\ner={\bf r}_1(t),t\in [0,2\pi)\}$$
in terms of a (global) positively oriented twice continuously differentiable $2\pi$-periodic parametrization ${\bf r}_1:[0,2\pi)\to\R^2$. More general piecewise smooth curves $\Gamma_1$ admitting local (patch-based) curve parametrizations as well as multiply connected obstacles~$\Omega$ can be easily incorporated in our approach.

Here,  $k_1=\omega\sqrt{\epsilon_1\mu_1}>0$ is the wavenumber of the exterior domain $D'=\R^2\setminus\overline D$ with  permittivity $\epsilon_1>0$ and  permeability $\mu_1>0$. Inside the penetrable array~$D$, on the other hand, the wavenumber is given by $k_2=\omega\sqrt{\epsilon_2\mu_2}$ in terms of $\mu_2>0$ and~$\epsilon_2$ which is allowed to be a complex number satisfying $\imag \epsilon_2\geq 0$.

\begin{figure}[h!]
\centering	
\includegraphics[scale=0.95]{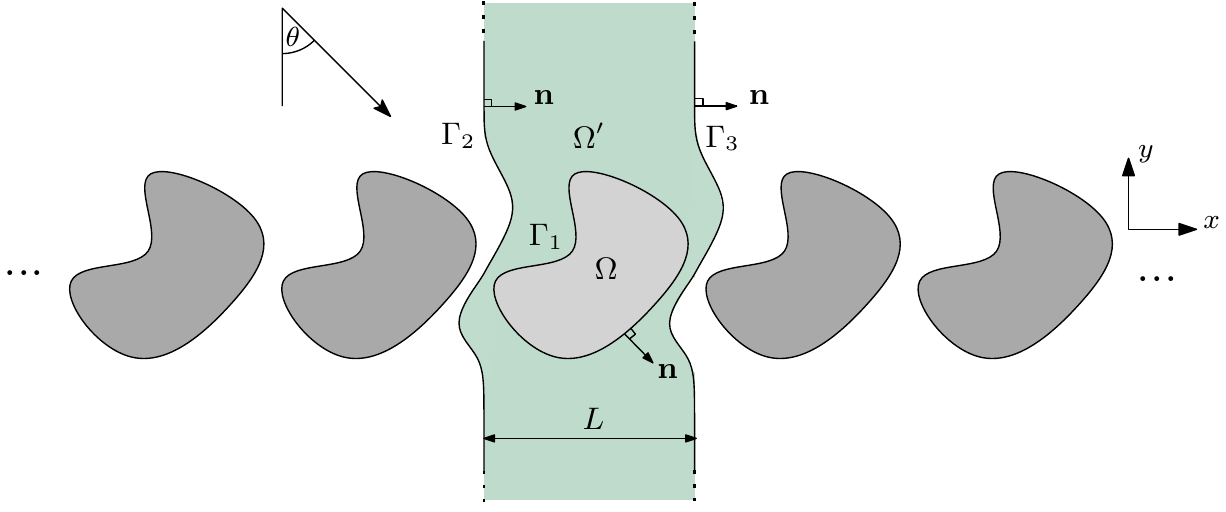}
\caption{Depiction of the quasi-periodic domain and the curves used in the derivation of the boundary integral equation formulation.}\label{fig:array}
\end{figure}

The sought total field $u:\R^2\to\C$, $u\in (C^2(D')\cap C^1(\overline{D'}))\cup (C^2(D)\cap C^1(\overline D))$, is the transverse component of the total electric field in TE polarization (resp. magnetic field in TM polarization) and solves 
\begin{subequations}\begin{equation}
\Delta u+k_1^2u = 0\quad\mbox{in}\quad D'\quad\mbox{and}\quad \Delta u+k_2^2u = 0\quad\mbox{in}\quad D.\label{eq:Helm_Eq}
\end{equation}
Additionally, it satisfies the quasi-periodicity condition
\begin{equation}
u(x+L,y)=\gamma u(x,y),\quad \gamma:=\e^{i\alpha L},\quad (x,y)\in \R^2,\label{eq:quasi-per}
\end{equation}
and the transmission conditions
\begin{equation}
 u|_{\p D}^+= u|_{\p D}^-\quad\mbox{and}\quad \p_n u|_{\p D}^+=\eta\p_n u|_{\p D}^-,\label{eq:trans_conditions}
\end{equation}\label{eq:trans_problem}\end{subequations}
where  $\eta:=\mu_{1}/\mu_2$ (resp. $\eta := \varepsilon_{1}/\varepsilon_2$). Here and in the sequel we use the notation
\begin{equation}\label{eq:traces}
u|_{\Gamma}^\pm(\ner)=\lim_{\delta\to 0^+}u(\ner\pm\delta {\nor}(\ner))\quad\text{and}\quad
\p_{n} u|_{\Gamma}^\pm(\ner)=\lim_{\delta\to 0^+} \nabla u(\ner\pm\delta \nor(\ner))\cdot \nor(\ner),\quad  \ner\in\Gamma,
\end{equation} 
 for the Dirichlet and Neumann traces, respectively, on a curve~$\Gamma$ with unit normal vector $\nor$. The precise orientation of $\nor$ is shown in Figure~\ref{fig:array}. 
 
As usual, we express the total field as
\begin{equation}
u=\begin{cases}
u^s+u^\inc&\text{in }\  D'\\
u^t&\text{in }\ D
\end{cases}\label{eq:decomp}
\end{equation}
in terms of the incident ($u^\inc$), transmitted ($u^t$), and scattered ($u^s$) fields, with the latter satisfying the Rayleigh expansion
\begin{subequations}\begin{equation}
u^s(x,y) = \left\{\begin{array}{ccl}\displaystyle\sum_{n\in\mathbb Z}B^+_n\e^{i(\alpha_n x+\beta_ny)}&\mbox{for}&  y>h^+:=\displaystyle\sup_{(x,y)\in\Omega}y\medskip\\
\displaystyle\sum_{n\in\mathbb Z}B^-_n\e^{i(\alpha_n x-\beta_ny)}&\mbox{for}&  y<h^-:=\displaystyle\inf_{(x,y)\in\Omega}y\end{array}\right.\label{eq:rad_cond_a}
\end{equation}
above ($y>h^+$) and below ($y<h^-$) the infinite array $D$, where  
\begin{equation}\label{eq:exp_coeff}
\alpha_n := \alpha +n\frac{2\pi}{L}\quad\mbox{and}\quad \beta_n:=\left\{\begin{array}{ccc}\sqrt{k_1^2-\alpha_n^2}&\mbox{if}&\alpha_n^2\leq k_1^2\medskip\\
i\sqrt{\alpha_n^2-k_1^2}&\mbox{if}&\alpha_n^2> k_1^2.\end{array}\right.
\end{equation}\label{eq:rad_cond}\end{subequations}

As it turns out it is convenient to distinguish the following three integer sets:
\begin{subequations}
\begin{align}
\mathcal{U}:=&\{n\in\Z: \alpha_n^2<k^2_1\}\\
\mathcal{V}:=&\{n\in\Z: \alpha_n^2>k^2_1\}\\
\mathcal{W}:=&\{n\in\Z: \alpha_n^2=k^2_1\}.
\end{align}\end{subequations}
According to our time convention,  it holds that for $n\in \mathcal U$  the modes 
\begin{equation}\label{eq:prog_modes}
u_n^+(x,y) := \e^{i\alpha_n x+i\beta_n y}\andtext u_n^{-}(x,y):= \e^{i\alpha_n x-i\beta_n y}
\end{equation} in~\eqref{eq:rad_cond_a} are upgoing and downgoing propagative planewaves, respectively. For $n\in\mathcal V$, in turn,  $u_n^+$ (resp. $u_n^{-}$) correspond to evanescent planewaves; they decay exponentially when $y\to +\infty$ (resp. $y \to -\infty$) while they grow exponentially as $y\to-\infty$ (resp. $y\to+\infty$). 

In turn, the set of integers 
$\mathcal{W}=\mathcal{W}(k_1,\alpha,L)=\left\{n \in \mathbb{Z}: \left(\alpha+2\pi n/L \right)^{2}=k_1^{2}\right\}=\left\{n \in \mathbb{Z}: \beta_{n}=0\right\}$, corresponds to the so-called \textit{Rayleigh-Wood anomaly configurations}~\cite{fernandez2020wave,bruno2017rapidly,bruno2020evaluation}.
For such $n$ values it holds that  
\begin{subequations}\begin{equation}
u_n(x,y):=u_n^+(x,y)=u_n^-(x,y)=\e^{i\alpha_n x}
\end{equation}
is a planewave that propagates parallel to the array along the $x$ axis. As it follows from separation of variables, the additional quasi-periodic homogeneous solution of the Helmholtz equation is given by the degenerated solution
\begin{equation}
v_n(x,y) := y\e^{i\alpha_n x},\quad n\in\mathcal W.
\label{eq:wood_mode}\end{equation}\label{eq:modes}\end{subequations}

Typically, the Rayleigh series~\eqref{eq:rad_cond} serves as the radiation condition for the quasi-periodic scattered field $u^s$. 
Alternatively, however,  such a radiation condition can be expressed in a less direct form by projecting the scattered field onto the non-radiative modes. This form of the radiation condition will turn out to be more suitable for our boundary integral equation formulation. To find it, we first note that since $u^s$ solves the homogeneous Helmholtz equation $\Delta u^s+k_1^2 u_s=0$ in $\Omega'$ and is quasi-periodic, it formally admits the general series expansion
\begin{equation}
u^s = \left\{\begin{array}{ccl}\displaystyle\sum_{n\in\mathcal U\cup\mathcal V}\left\{B^+_nu^{+}_n+C^+_nu^{-}_n\right\}+\sum_{n\in\mathcal W}\left\{B_n^+u_n+C_n^+v_n\right\}&\mbox{for}&  y>h^+\medskip\\
\displaystyle\sum_{n\in\mathcal U\cup\mathcal V}\left\{B^-_nu_n^{-}+C^-_nu_n^{+}\right\}+\sum_{n\in\mathcal W}\left\{B_n^-u_n+C_n^-v_n\right\}&\mbox{for}&  y<h^-.\end{array}\right.\label{eq:general}
\end{equation} The fact that $u^s$ is radiative and bounded, in the sense of~\eqref{eq:rad_cond}, then implies that $C^\pm_n=0$ for all $n\in\mathbb Z$. Therefore, computing these coefficients by projecting $u^s(\cdot,\pm h)$ and $\p_y u^s(\cdot,\pm h)$for $h>\max\{h^+,-h^-\}$ onto $\e^{i\alpha_n\cdot}$, we conclude that the radiation condition~\eqref{eq:rad_cond} can be equivalently enforced by requesting $u^s$ to satisfy~\eqref{eq:general} and
\begin{equation}\label{eq:rad_cond_2}\begin{split}
\frac{1}{L}\int_{-\frac{L}2}^{\frac{L}2}\left\{\p_y u^s(x,\pm h)\mp i\beta_n u^s(x,\pm h)\right\}\e^{-i\alpha_nx}\de x&=0\\
&=\begin{cases}\mp 2i\beta_n \e^{-i\beta_n h}C_n^\pm&\text{if}\ n\in\mathcal U\cup\mathcal V\\C_n^\pm&\text{if}\ n\in\mathcal W.\end{cases}\end{split}
\end{equation}

Unlike most of previous works, our direct BIE formulation, which  is presented in Section~\ref{sec:BIE}, leverages the representation formula~\eqref{eq:GreenFormula} of the scattered field $u^s$ 
that uses the standard free-space Helmholtz Green function $G_{k_1}$ given by 
\begin{equation}
G_{k}(\ner,\ner') := \frac{i}{4}H_0^{(1)}(k|\ner-\ner'|)\quad(\ner\neq\ner').
\label{eq:free-space}\end{equation} 
This formula, which is derived in the next section, however,  involves integration over the boundary $\Gamma_2\cup\Gamma_3$ of the unit cell domain
\begin{equation}\label{eq:unit_cell}
 R=\{(x,y)\in\R^2: y={\mathsf y}_{2}(t), {\mathsf x}_{2}(t)<x<{\mathsf x}_{2}(t)+L,t\in\R\}.
\end{equation}
where the infinite parallel curves
\begin{equation}\label{eq:par_curves}
\Gamma_2 := \{\ner\in\R^2: \ner={\bf r}_2(t),t\in \R\}\quad\text{and}\quad \Gamma_3:=\{\ner\in\R^2: \ner={\bf r}_{2}(t)+L{\bf e}_1,t\in \R\}
\end{equation}
 are parametrized by the smooth function  ${\bf r}_2(t)=({\mathsf x}_{2}(t),{\mathsf y}_{2}(t))$.
These curves are assumed  to extend infinitely along the $y$-axis while satisfying the condition  $ {\mathsf y}_{2}(t)= t $ for $|t|>\max\{h^+,-h^-\}$ and not intercepting~$\Gamma_1$ (see Figure~\ref{fig:array}). These  conditions simplify the analysis of the windowed integral operators in Sections~\ref{sec:WGF}-\ref{sec:corrected_BIE} and Appendix~\ref{rem:sub_al_decay}.

\section{Green's representation formulae} \label{sec:green_representation}
We start off this section by deriving integral representation formulae for the total field $u$ in $\Omega$ and the scattered field $u^s$ in $\Omega':=R\setminus\overline\Omega$, where $R$ is unit cell-domain~\eqref{eq:unit_cell},  using the free-space Helmholtz Green function~\eqref{eq:free-space}.

\begin{figure}[h!]
\centering	
\includegraphics[scale=0.95]{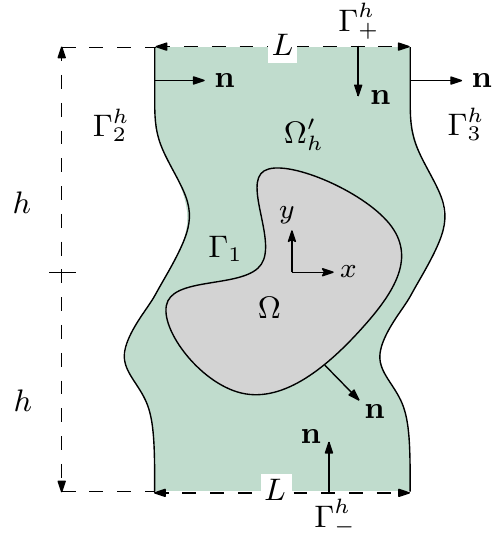}
\caption{Depiction of the curves involved in the derivation of Green's representation formula.}\label{fig:green}
\end{figure}

Let us consider the bounded domain $\Omega^\prime_h=R_h\cap\Omega'$ where $R_h = \{(x,y)\in R: |y|<h\}$ and $h>\max\{h^+,-h^-\}$; see Figure~\ref{fig:green} and~\eqref{eq:rad_cond_a} for the definition of $h^\pm$. Applying Green's third identity we have that for any fixed target point $\ner=\left(x,y\right) \in \Omega^\prime_h$ it holds that
\begin{equation*}\label{eq:green_start}\begin{split}
\lf(\int\limits_{\Gamma_1}+\int\limits_{\Gamma^h_+}+\int\limits_{\Gamma^h_-}+\int\limits_{\Gamma^h_2}-\int\limits_{\Gamma^h_3}\rg) \left\{u^s(\ner') \frac{\p G_{k_1}(\ner,\ner')}{\p n(\ner')} -\p_n u^s(\ner')G_{k_1}(\ner,\ner') \right\}\de s(\ner')=\\\begin{cases}u^{s}(\ner)&\text{if }\ner\in\Omega'_h\\0&\text{if }\ner\in\Omega\end{cases}
\end{split}\end{equation*}
where integration is carried out over the multiply connected curve $\p\Omega'_h$ that comprises $\Gamma_1$, $\Gamma_2^h$, $\Gamma_3^h$, and the horizontal lines 
$\Gamma^h_\pm=\{(x,y)\in R: y=\pm h\}$  with normals pointing toward the interior of $\Omega'_h$.

Applying the Cauchy-Schwartz inequality, we have 
\begin{equation*}\begin{split}
\lf|\int_{\Gamma_\pm^h} u^s(\ner') \frac{\p G_{k_1}(\ner,\ner')}{\p n(\ner')}  \de s(\ner')\rg|&\leq \lf(\int_{\Gamma_\pm^h}|u^s|^2\de s\rg)^{1/2}\lf(\int_{\Gamma_\pm^h}\lf|\frac{\p G_{k_1}(\ner,\ner')}{\p n(\ner')}\rg|^2\de s(\ner')\rg)^{1/2},\\
\lf| \int_{\Gamma_\pm^h} \p_{n} u^s(\ner')G_{k_1}(\ner,\ner') \de s(\ner')\rg|&\leq\lf(\int_{\Gamma_\pm^h}|\p_n u^s|^2\de s\rg)^{1/2}\!\lf(\int_{\Gamma_\pm^h}\lf| G_{k_1}(\ner,\ner')\rg|^2\de s(\ner')\rg)^{1/2}.
\end{split}\end{equation*}
Therefore,  from the uniform boundedness of $u^s$ and $\p_n u^s$ on $\Gamma_\pm^h$ for all $h>\max\{h^+,-h^-\}$ (which follows from~\eqref{eq:rad_cond}), and the large-argument asymptotic expansion of the Hankel functions~(see e.g.~\cite{Abramowitz1966Handbook}), we obtain 
\begin{equation*}\begin{split}
\lf|\int_{\Gamma_\pm^h} u^s(\ner') \frac{\p G_{k_1}(\ner,\ner')}{\p n(\ner')}  \de s(\ner')\rg|\lesssim&~ \frac{1}{\sqrt h}\to 0\ \text{ and }\\
\ \lf| \int_{\Gamma_\pm^h} \p_{n} u^s(\ner')G_{k_1}(\ner,\ner') \de s(\ner')\rg|\lesssim&~\frac{1}{\sqrt h}\to 0\quad\text{as}\quad  h\to \infty.
\end{split}\end{equation*}

Taking then the limit as $h\to\infty$ in the remaining integrals over $\Gamma_2^h$ and $\Gamma_3^h$ we arrive at the integral representation formula
\begin{equation}\label{eq:GreenFormula}\begin{split}
\lf(\ \int\limits_{\Gamma_1}+\int\limits_{\Gamma_2}-\int\limits_{\Gamma_3}\rg) \left\{u^s(\ner') \frac{\p G_{k_1}(\ner,\ner')}{\p n(\ner')} -\p_n u^s(\ner')G_{k_1}(\ner,\ner') \right\}\de s(\ner')=\\
\begin{cases}u^{s}(\ner)&\text{if }\ner\in\Omega'\\0&\text{if }\ner\in\Omega.\end{cases}
\end{split}\end{equation}

Note that only the boundedness of the scattered field and its gradient was used in these derivations, not the radiation condition. More, precisely, any bounded homogenous solutions of the Helmholtz equation in $\Omega'$ admits the integral representation~\eqref{eq:GreenFormula}. In particular, it can be shown (cf.~\cite{DeSanto:1998tq}) that upgoing and downgoing planewaves~\eqref{eq:prog_modes} for $n\in\mathcal U$ as well as horizontally propagating modes~\eqref{eq:wood_mode} for $n\in\mathcal W$ satisfy
\begin{equation}\label{eq:GreenFormula2}
\lf(\ \int\limits_{\Gamma_2}-\int\limits_{\Gamma_3}\rg) \left\{u^{\pm}_n(\ner') \frac{\p G_{k_1}(\ner,\ner')}{\p n(\ner')} -\p_n u^{\pm}_n(\ner')G_{k_1}(\ner,\ner') \right\}\de s(\ner')=u_n^{\pm}(\ner),\ \ \ner\in R.
\end{equation}

Finally, applying the standard Green's third identity inside $\Omega$, we readily obtain the representation formulae
\begin{equation}\label{eq:GreenFormulaInc}
-\int_{\Gamma_1} \left\{u^\inc(\ner') \frac{\p G_{k_1}(\ner,\ner')}{\p n(\ner')} -\p_n u^\inc(\ner')G_{k_1}(\ner,\ner') \right\}\de s(\ner')=\begin{cases}0&\text{if }\ner\in\Omega'\\u^\inc(\ner)&\text{if }\ner\in\Omega\end{cases}
\end{equation}
and 
\begin{equation}\label{eq:GreenFormulaTrans}
-\int_{\Gamma_1} \left\{u^t(\ner') \frac{\p G_{k_2}(\ner,\ner')}{\p n(\ner')} -\p_n u^t(\ner')G_{k_2}(\ner,\ner') \right\}\de s(\ner')=\begin{cases}0&\text{if }\ner\in\Omega'\\u^t(\ner)&\text{if }\ner\in\Omega.\end{cases}
\end{equation}
for the incident and transmitted fields inside $\Omega$.

\section{Parametrized integral operators}\label{sec:param_op}

This section presents the notation and the main properties of the Helmholtz layer potentials and boundary integral operators used in the construction of the direct boundary integral equation in Section~\ref{sec:BIE}.

For a sufficiently regular density function $\varphi:\Gamma_i\to\C$, with $\Gamma_i$ being one of the curves introduced above, we define the single- and double-layer potentials as
\begin{equation}\begin{aligned}
(S\!L_j^{i}\varphi)(\ner) :=&~ \int_{\Gamma_i} G_{k_j}(\ner,\ner')\varphi(\ner')\de s(\ner')\quad\mbox{and}\\
(D\!L_j^{i}\varphi)(\ner) :=&~ \int_{\Gamma_i} \frac{\p G_{k_j}(\ner,\ner')}{\p n({\ner}')}\varphi(\ner')\de s(\ner'),\quad \ner\in\R^2\setminus\Gamma,
\end{aligned}\label{eq:layer_potentials}\end{equation}
respectively, corresponding to the wavenumbers $k_j$, $j=1,2.$ For the potentials associated with the unbounded curves $\Gamma_i$,  $i=2,3$, the boundary integrals in~\eqref{eq:layer_potentials} should be interpreted as improper conditionally convergent integrals. 

As is well-known (cf. \cite{COLTON:1983}), these potentials satisfy the jump relations
\begin{equation}
\begin{aligned}
  \left.S\!L_j^{i}\varphi\right|_{\Gamma_\ell}^\pm =&~V_j^{\ell,i}\varphi,&\qquad \left.\p_{n} (S\!L_j^{i}\varphi)\right|_{\Gamma_\ell}^\pm =
\displaystyle \mp\delta_{i,\ell}\frac{\varphi}{2} + \widetilde K_j^{\ell,i}\varphi,
\\
\left.\p_{n}(D\!L_j^{i}\varphi)\right|_{\Gamma_\ell}^\pm =&~ W_j^{\ell,i}\varphi,&\qquad \left.D\!L_j^{i}\varphi\right|_{\Gamma_\ell}^\pm = 
\displaystyle \pm\delta_{i,\ell}\frac{\varphi}{2} + K_j^{\ell,i}\varphi,
 \end{aligned}\label{eq:jump_conds}
\end{equation}
(see the notation introduced in~\eqref{eq:traces}) where $\delta_{ij}$ denotes the Kronecker delta, and  $V_j^{\ell,i}$, $K_j^{\ell,i}$, $\widetilde K_j^{\ell,i}$, $W_j^{\ell,i}$ above are, respectively, the single-layer, double-layer, adjoint double-layer, and hypersingular operators with wavenumber $k_j$, target surface $\Gamma_\ell$, and source surface $\Gamma_i$;  these operators are formally defined as
\begin{subequations}\begin{align}
 \big(V_j^{\ell, i}\varphi\big)(\ner) :=&~ \int_{\Gamma_i} G_{k_j}(\ner,\ner')\varphi(\ner')\de s(\ner') \\
  \big(K_j^{\ell, i}\varphi\big)(\ner) :=&~ \int_{\Gamma_i} \frac{\p G_{k_j}(\ner,\ner')}{\p n(\ner')}\varphi(\ner')\de s(\ner')\medskip\\
 \big(\widetilde K_j^{\ell, i}\varphi\big)(\ner) :=&~ \int_{\Gamma_i} \frac{\p G_{k_j}(\ner,\ner')}{\p n(\ner)}\varphi(\ner')\de s(\ner') \\
 \big(W_j^{\ell, i}\varphi\big) (\ner):=&~ {\rm f.p.}\!\!\int_{\Gamma_i} \frac{\p^2 G_{k_j}(\ner,\ner')}{\p n(\ner)\p n(\ner')}\varphi(\ner')\de s(\ner'),\quad \ner\in\Gamma_\ell \label{eq:hyp_op}
\end{align}\label{eq:operators}\end{subequations}
 (the letters f.p. indicates that the integral in the definition of the hypersingular operator~\eqref{eq:hyp_op} has to be interpreted as a Hadamard finite part integral).  

In order to deal with the quasi-periodicity condition, it will be convenient to work with the layer potentials~\eqref{eq:layer_potentials} and the integral operators~\eqref{eq:operators} in parametric form. For a sufficiently smooth density function $\varphi:\Gamma_i\to \C$, we let $\phi = \varphi\circ {\bf r}_i:I_i\to\C$, with $I_1=[0,2\pi)$ and $I_{2}=I_3=\R,$ and define the parametric layer potentials as $(\mathcal S_j^{i}\phi)(\ner) := (S\!L_j^i\varphi)(\ner)$ and  $(\mathcal D_j^{i}\phi)(\ner) :=(D\!L_j^i\varphi)(\ner)$, or, more explicitly as
\begin{equation}\begin{aligned}
(\mathcal {S}_j^{i}\phi)(\ner) =&~\frac{i}{4}\int_{I_i}H_0^{(1)}(k_j|\ner-{\bf r}_i(\tau)|)\phi(\tau)|{\bf r}_i'(\tau)|\de\tau,\\
(\mathcal D_j^{i}\phi)(\ner) =&~\frac{ik_j}{4}\int_{I_i}H_1^{(1)}(k_j|\ner-{\bf r}_i(\tau)|)\frac{(\ner-{\bf r}_i(\tau))\cdot{\bf n}_i(\tau)}{|\ner-{\bf r}_i(\tau)|}\phi(\tau)|{\bf r}_i'(\tau)|\de\tau,  \  \ \ner\in\R^2\setminus\Gamma
\end{aligned}\label{eq:param_pots}\end{equation}
where ${\bf n}_i = (\mathsf y_i',-\mathsf x'_i)/|{\bf r}'_i|$ denotes the parametrized unit normal vector to the curve~$\Gamma_i$.

Similarly, the parametric boundary integral operators are defined as ${\mathsf{V}}_j^{\ell,i}\phi=(V_j^{\ell,i}\varphi)\circ {{\bf r}}_\ell$,  ${\mathsf K}_j^{\ell,i}\phi=(K_j^{\ell,i}\varphi)\circ {{\bf r}}_\ell$, $\widetilde{\mathsf K}_j^{\ell,i}\phi=(\widetilde K_j^{\ell,i}\varphi)\circ {{\bf r}}_\ell$, and  ${\mathsf W}_j^{\ell,i}\phi=(W_j^{\ell,i}\varphi)\circ {{\bf r}}_\ell$. For self-containess we write them in extensive as
\begin{equation}\begin{aligned}
({\mathsf{V}}_j^{\ell,i}\phi)(t) =&\int_{I_i}Q_{V,k}^{\ell,i}(t,\tau)\phi(\tau)|{\bf r}_i'(\tau)|\de\tau, &\qquad({\mathsf{K}}_j^{\ell,i}\phi)(t) =&\int_{I_i}Q_{K,k}^{\ell,i}(t,\tau)\phi(\tau)|{\bf r}_i'(\tau)|\de\tau, \\
(\widetilde{\mathsf K}_j^{\ell,i}\phi)(t) =&\int_{I_i}Q_{\widetilde K,k}^{\ell,i}(t,\tau)\phi(\tau)|{\bf r}_i'(\tau)|\de\tau, &\qquad (\mathsf W_j^{\ell,i}\phi)(t)=&~{\rm f.p.}\!\!\int_{I_i}Q_{W,k}^{\ell,i}(t,\tau)\phi(\tau)|{\bf r}_i'(\tau)|\de\tau
\end{aligned}\label{eq:op_param}\end{equation}
where letting $\bol R_{\ell,i} = {\bf r}_\ell(t)-{\bf r}_i(\tau)$ and $R_{\ell,i} = |{\bf r}_\ell(t)-{\bf r}_i(\tau)|$ the integral kernels above are given by
\begin{subequations}\begin{align}
Q^{\ell,i}_{V,j}(t,\tau):=&~\frac{i}{4}H_0^{(1)}(k_jR_{\ell,i})\\
Q^{\ell,i}_{K,j}(t,\tau):=&~\frac{ik_j}{4}H_1^{(1)}(k_jR_{\ell,i})\frac{\bol R_{\ell,i}\cdot{\bf n}_i(\tau)}{R_{\ell,i}}\\
Q^{\ell,i}_{\widetilde K,j}(t,\tau):=&-\frac{ik_j}{4}H_1^{(1)}(k_jR_{\ell,i})\frac{\bol R_{\ell,i}\cdot{\bf n}_\ell(t)}{R_{\ell,i}}\\
 Q^{\ell,i}_{W,j}(t,\tau):=&~\frac{i k_j}{4}\lf(\frac{H_{1}^{(1)}(k_j R_{\ell,i})}{R_{\ell,i}} \mathbf{n}_\ell(t) \cdot \mathbf{n}_i(\tau)\rg.+\label{eq:param_ker_hyper}\\
&\lf.\left\{k_j R_{\ell,i} H_{0}^{(1)}(k_j R_{\ell,i})-2 H_{1}^{(1)}(k_j R_{\ell,i})\right\} \frac{\bol R_{\ell,i} \cdot \mathbf{n}_{i}(\tau)\bol R_{\ell,i} \cdot \mathbf{n}_\ell(t)}{R_{\ell,i}^{3}}\rg)\nonumber
\end{align}\label{eq:param_ker}\end{subequations}
for $i,\ell = 1,2,3$ and $j=1,2$.

The following simple result greatly simplifies the final form of the our boundary integral equation system:
\begin{proposition}\label{prop:op_indeitities}
The  identities 
\begin{equation}
\mathsf{V}_{1}^{2,2} =\mathsf V_1^{3,3},\quad \mathsf K_{1}^{2,2} = \mathsf K_1^{3,3},\quad \mathsf{\widetilde K}_{1}^{2,2} = \mathsf{\widetilde K}_1^{3,3}\andtext\mathsf W_{1}^{2,2} = \mathsf W_1^{3,3}\label{eq:op_id1}\end{equation} hold for the parametrized integral operators defined in~\eqref{eq:op_param} associated with the parallel curves $\Gamma_2$ and $\Gamma_3$ defined in~\eqref{eq:par_curves}. 
Furthermore, 
\begin{equation}
\mathsf{V}_{1}^{2,3} =\mathsf V_1^{3,2},\quad \mathsf  K_{1}^{2,3} = -\mathsf K_1^{3,2},\quad  \mathsf{\widetilde K}_{1}^{2,3} =  -\mathsf{\widetilde K}_1^{3,2}\andtext  \mathsf  W_{1}^{2,3} =\mathsf  W_1^{3,2}\label{eq:op_id2}\end{equation} 
in the particular case when $\Gamma_2$ and $\Gamma_3$ are straight vertical lines with constant unit normal ${\bf e}_1$.
\end{proposition}

\begin{proof}
The first part~\eqref{eq:op_id1} follows directly from the fact ${\bf r}_3(t)={\bf r}_2(t)+L{\bf e}_1$, $t\in\R$ and hence $R_{2,2}=R_{3,3}$, $\bol R_{2,2}=\bol R_{3,3}$ and ${\bf n}_2={\bf n}_3$ in~\eqref{eq:param_ker}. For the second part~\eqref{eq:op_id2} the proof follows  from writing the integral kernels~\eqref{eq:param_ker} using the parametrizations ${\bf r}_2(t) = -\frac{L}2{\bf e}_1 +\mathsf y_2(t){\bf e}_2$ for $\Gamma_2$ and ${\bf r}_3(t) = \frac{L}2{\bf e}_1+\mathsf y_2(t){\bf e}_2$ for $\Gamma_3$ that have a constant unit normal ${\bf n}_2(t)={\bf n}_3(t)=
{\bf e}_1$ for all $t\in\R$. Therefore,  in view of the fact that $R_{2,3} = \sqrt{L^2+(\mathsf y_2(t)-\mathsf y_2(\tau))^2}=R_{3,2}$, $\bol R_{2,3}\cdot {\bf e}_1 =-L= -\bol R_{3,2}\cdot {\bf e}_1$ and $|{\bf r}'_3|=|{\bf r}'_3|$ in this case, the identities in~\eqref{eq:op_id2} readily follow.\end{proof}

\section{Boundary integral equation formulation}\label{sec:BIE}

A direct boundary integral equation formulation for the quasi-periodic transmission problem presented in Section~\ref{sec:Inte_Eq_For} is derived in this section. Our strategy lies in recasting the problem as a (formally) second-kind system of boundary integral equations for the interior traces of the total field on $\Gamma_1$ and for the traces of the scattered field on the unbounded curves $\Gamma_2$ and $\Gamma_3$. We follow here the Kress-Roach approach~\cite{kress1978transmission} (also known as M\"uller's formulation~\cite{muller2013foundations} for its 3D electromagnetic version) which yields two second-kind integral equations from enforcing the transmission conditions~\eqref{eq:trans_conditions} on $\Gamma_1$. The remaining two equations, on the other hand, are derived from the representation formula~\eqref{eq:GreenFormula} that is used to suitably combine the traces of the scattered field on $\Gamma_2$ and $\Gamma_3$, to obtain second-kind equations that account for the quasi-periodicity of the scattered field. One salient advantage of our approach is that the resulting integral operators are expressed in terms weakly-singular and smooth kernels that can integrated with high precision using global trigonometric quadrature rules.

We start off by noting that by virtue of the quasi-periodicity condition~\eqref{eq:quasi-per} the traces $u^s|_{\Gamma_3}^-$ and $\p_n u^s|_{\Gamma_3}^-$ can be expressed in terms of $u^s|_{\Gamma_2}^+$ and $\p_n u^s|_{\Gamma_2}^+$. Indeed, using the curve parametrizations ${\bf r}_i:I_i\to\Gamma_i$ for the curves $\Gamma_i$, $i=1,2,3$, we  have the scattered field traces on the $\Gamma_2$  and $\Gamma_3$ satisfy
\begin{equation}\label{eq:trace_relations}
u^s|_{\Gamma_3}^-\circ{\bf r}_3 = \gamma u^s|_{\Gamma_2}^+\circ{\bf r}_2 \quad\mbox{and}\quad \p_{n} u^s |_{\Gamma_3}^- \circ {\bf r}_3=\gamma \p_n u^s|_{\Gamma_2}^+\circ{\bf r}_2\quad (\gamma=\e^{i\alpha L}),
\end{equation} 
where the facts that ${\bf r}_3={\bf r}_2+L{\bf e}_1$ and that the curves thus parametrized share the same unit normal ${\bf n}_2={\bf n}_3= (\mathsf y_2',-\mathsf x'_2)/|{\bf r}'_2|$ were used. It hence follows from~\eqref{eq:trace_relations} that only the parametrized traces
\begin{equation}\label{eq:unknowns}\begin{aligned}
\phi_1=&~u^t|_{\Gamma_1}^-\circ{\bf r}_1:[0,2\pi)\to\C\qquad\qquad& \phi_2=&~\p_n u^t|_{\Gamma_1}^-\circ{\bf r}_1:[0,2\pi)\to\C\\
\phi_3=&~u^s|_{\Gamma_2}^+\circ{\bf r}_2:\R\to\C\qquad\qquad& \phi_4=&~\p_{n} u^s|_{\Gamma_2}^+\circ{\bf r}_2:\R\to\C
\end{aligned}\end{equation}
are needed in order to retrieve the fields by means of the representation formulae~\eqref{eq:GreenFormula},~\eqref{eq:GreenFormulaTrans} and~\eqref{eq:GreenFormulaInc}. 

Indeed, since by the transmission conditions~\eqref{eq:trans_conditions} we have $u^s|_{\Gamma_1}^+ = u^t|_{\Gamma_1}^--u^\inc|_{\Gamma_1} $ and $\p_n u^s|_{\Gamma_1}^+ = \eta\p_n u^t|_{\Gamma_1}^--\p_n u^\inc|_{\Gamma_1} $, the representation formulae~\eqref{eq:GreenFormula} and~\eqref{eq:GreenFormulaInc} can be combined with~\eqref{eq:trace_relations} to obtain the following integral representation of the scattered field
\begin{equation}\begin{split}
u^s(\ner) = \big(\mathcal D_1^{1}\phi_1\big)(\ner)  -\eta\big(\mathcal S_1^{1}\phi_2\big)(\ner)+ \big(\mathcal D_1^{2}\phi_3\big)(\ner)  -\big(\mathcal S_1^{2}\phi_4\big)(\ner)\\
-\gamma\lf\{\big(\mathcal D_1^{3}\phi_3\big)(\ner) -\big(\mathcal S_1^{3}\phi_4\big)(\ner)\rg\},& \quad \ner\in\Omega' \label{eq:GF_scat}
\end{split}\end{equation}
where we used the parametrized form of the layer potentials~\eqref{eq:param_pots}. Similarly, the transmitted field in~\eqref{eq:GreenFormulaTrans} can be expressed as
\begin{equation}
u^t(\ner) = -\big(\mathcal D_2^{1}\phi_1\big)(\ner) +\big(\mathcal S_2^{1}\phi_2\big)(\ner), \quad \ner\in\Omega \label{eq:GF_trans}
\end{equation} in terms of the unknown densities~\eqref{eq:unknowns}.

Letting 
$$
f = u^{\inc}\circ{\bf r}_1\andtext g ={\bf n}_1\cdot \nabla u^{\inc}\circ{\bf r}_1
$$
we have that a direct application of the jump conditions~\eqref{eq:jump_conds} to evaluate~\eqref{eq:GF_scat} and its normal derivative on $\Gamma_1$,  yields the equations
\begin{subequations}\begin{eqnarray}
-f+\frac{\phi_1}{2}  &=&\mathsf K^{1,1}_1\phi_1- \eta \mathsf V^{1,1}_1\phi_2+(\mathsf K^{1,2}_1-\gamma \mathsf K_1^{1,3})\phi_3- (\mathsf V^{1,2}_1-\gamma \mathsf V^{1,3}_1)\phi_4\label{eq:B1}\\
-g +\frac{\eta}2\phi_2 &=&\mathsf W^{1,1}_1\phi_1- \eta\widetilde{\mathsf K}^{1,1}_1\phi_2+(\mathsf W^{1,2}_1-\gamma \mathsf W^{1,3}_1)\phi_3- (\widetilde{\mathsf K}^{1,2}_1-\gamma \widetilde{\mathsf K}^{1,3}_1)\phi_4\label{eq:B2}
\end{eqnarray}\label{eq:gamma_1}\end{subequations}
which hold in $[0,2\pi)$. Similarly, using~\eqref{eq:jump_conds} to evaluate the transmitted field~\eqref{eq:GF_trans} and its normal derivative on $\Gamma_1$, we obtain 
\begin{subequations}\begin{align}
 \frac{\phi_1}{2} =& -\mathsf K^{1,1}_2\phi_1+ \mathsf V^{1,1}_2\phi_2\label{eq:A1}\\
 \frac{\phi_2}{2} =& -\mathsf W^{1,1}_2\phi_1+ \widetilde{\mathsf K}^{1,1}_2\phi_2\label{eq:A2}
\end{align}\label{eq:trans_field}\end{subequations} which hold in $[0,2\pi)$.

Therefore, adding~\eqref{eq:A1} to~\eqref{eq:B1} and adding~\eqref{eq:A2} to~\eqref{eq:B2} we arrive at the following integral equations 
\begin{equation}\label{eq:eq_set_1}
\phi_1+\sum_{q=1}^4{\mathsf M}_{1,q}\phi_q=f \andtext\lf(\frac{1+\eta}{2}\rg)\phi_2 +\sum_{q=1}^4{\mathsf M}_{2,q}\phi_q=g\quad\text{in}\quad [0,2\pi)
\end{equation}
where 
\begin{equation}\begin{array}{llll}
\mathsf M_{1,1}:= \mathsf K^{1,1}_2-\mathsf K^{1,1}_1,& \mathsf M_{1,2} :=  \eta\mathsf V^{1,1}_1-\mathsf V^{1,1}_2,&\mathsf M_{1,3}:=\gamma\mathsf K^{1,3}_1-\mathsf K^{1,2}_1,& \mathsf M_{1,4}:=\mathsf V^{1,2}_1-\gamma\mathsf V_1^{1,3},\smallskip\\
\mathsf M_{2,1} :=\mathsf W^{1,1}_2-\mathsf W^{1,1}_1,&  \mathsf M_{2,2}:=\eta\mathsf{\widetilde K}^{1,1}_1-\mathsf{\widetilde K}^{1,1}_2,&\mathsf M_{2,3}:=\gamma\mathsf W^{1,3}_1-\mathsf W^{1,2}_1,&\mathsf M_{2,4}:=\mathsf{\widetilde K}^{1,2}_1-\gamma \mathsf{\widetilde K}^{1,3}_1.
\end{array}\label{eq:blocks1}\end{equation}

\begin{remark}As mentioned above, all the integral operators in~\eqref{eq:blocks1} are weakly singular. Indeed, for instance,  the seemingly hypersingular operator $\mathsf M_{2,1}=\mathsf W^{1,1}_2-\mathsf W^{1,1}_1$ is weakly singular by virtue of the fact that hypersingular parametric kernel, defined in~\eqref{eq:param_ker_hyper}, can be expressed as
$$
Q^{1,1}_{W,j}(t,\tau)= \frac{\mathbf{n}_1(t) \cdot \mathbf{n}_1(\tau)}{2\pi R^{2}_{1,1}} +a_{j}(t, \tau) \log (|t-\tau|)+b_{j}(t, \tau), \quad t,\tau\in[0,2\pi),
$$
where $a_j,b_j: [0,2\pi)^2\to \C$ are smooth $2\pi$-periodic functions in both arguments~\cite{Kress:1995}.  Therefore, since the hypersingular static terms $\frac{\mathbf{n}_1(t) \cdot \mathbf{n}_1(\tau)}{2\pi R^{2}_{1,1}} $ cancels when we take the difference $Q^{1,1}_{W,2}-Q^{1,1}_{W,1}$,  the integral kernel of $\mathsf M_{2,1}$ features only a logarithmic singularity as $t\to\tau$.\qed \end{remark}

In order to find the two additional integral equations, we take the Dirichlet and Neumann traces~\eqref{eq:traces} of~\eqref{eq:GF_scat} on~$\Gamma_2$ and $\Gamma_3$ using the jump relations~\eqref{eq:jump_conds}, to obtain 
\begin{subequations}\begin{eqnarray}
\frac{\phi_3}{2}  &=&\mathsf K^{2,1}_1\phi_1- \eta \mathsf V^{2,1}_1\phi_2+( \mathsf K^{2,2}_1-\gamma  \mathsf K_1^{2,3})\phi_3- (\mathsf V^{2,2}_1-\gamma \mathsf V^{2,3}_1)\phi_4\label{eq:C1}\\
\frac{\phi_4}2 &=&\mathsf W^{2,1}_1\phi_1- \eta\widetilde{\mathsf K}^{2,1}_1\phi_2+(\mathsf W^{2,2}_1-\gamma \mathsf W^{2,3}_1)\phi_3- (\widetilde{\mathsf K}^{2,2}_1-\gamma \widetilde{\mathsf K}^{2,3}_1)\phi_4\label{eq:C2}\\
\gamma\frac{\phi_3}{2}  &=&\mathsf K^{3,1}_1\phi_1- \eta \mathsf V^{3,1}_1\phi_2+(\mathsf K^{3,2}_1-\gamma \mathsf K_1^{3,3})\phi_3- (\mathsf V^{3,2}_1-\gamma\mathsf V^{3,3}_1)\phi_4\label{eq:D1}\\
\gamma\frac{\phi_4}2 &=&\mathsf W^{3,1}_1\phi_1- \eta\widetilde{\mathsf K}^{3,1}_1\phi_2+(\mathsf W^{3,2}_1-\gamma \mathsf W^{3,3}_1)\phi_3- (\widetilde{\mathsf K}^{3,2}_1-\gamma \widetilde{\mathsf K}^{3,3}_1)\phi_4\label{eq:D2}
\end{eqnarray}\label{eq:gamma_3}\end{subequations}
which hold in $\R$. We then  combine these equations  to cancel  all the weakly-singular  ($\mathsf V_1^{i,i}, \mathsf K_1^{i,i}$, and $\mathsf {\widetilde K}_1^{i,i}$, $i=2,3$) and hypersingular ($\mathsf W_1^{i,i}$, $i=2,3$) operators. In detail,  multiplying \eqref{eq:C1} by $\gamma$ and adding it to~\eqref{eq:D1}, and multiplying \eqref{eq:C2} by $\gamma$ and adding it to~\eqref{eq:D2}, while using the identities in~\eqref{eq:op_id1}, we arrive at
\begin{equation}\label{eq:eq_set_2}
\gamma\phi_3+\sum_{q=1}^4 \mathsf M_{3,q}\phi_q=0\andtext \gamma\phi_4+\sum_{q=1}^4\mathsf M_{4,q}\phi_q=0\quad\text{in}\quad \R
\end{equation}
where
\begin{equation}\hspace{-0.2cm}\begin{array}{llll}
{\mathsf M}_{3,1} :=  -\gamma \mathsf K^{2,1}_1- \mathsf K^{3,1}_1,&\!\! \mathsf M_{3,2}= \eta (\gamma \mathsf V^{2,1}_1+\mathsf V^{3,1}_1),&\!\! \mathsf M_{3,3}:= \gamma^2\mathsf K_1^{2,3}-\mathsf K^{3,2}_1,&\!\!\mathsf M_{3,4}=  \mathsf V^{3,2}_1-\gamma^2\mathsf V^{2,3}_1,\smallskip\\
{\mathsf M}_{4,1} :=  -\gamma  \mathsf W^{2,1}_1- \mathsf W^{3,1}_1,&\!\!\mathsf M_{4,2}= \eta (\gamma \mathsf{\widetilde K}^{2,1}_1+\mathsf{\widetilde K}^{3,1}_1),&\!\! \mathsf M_{4,3}:= \gamma^2 \mathsf W^{2,3}_1-\mathsf W^{3,2}_1,&\!\!\mathsf M_{4,4}= \mathsf{\widetilde K}^{3,2}_1-\gamma^2 \mathsf{\widetilde K}^{2,3}_1.
\end{array}\label{eq:blocks2}\end{equation}

Clearly, the operators~\eqref{eq:blocks2} have smooth kernels, by virtue of the fact that integration and evaluation are carried out over different well-separated curves.

Finally, lumping the unknown density functions~\eqref{eq:unknowns} in a single vector $\bol\phi=\lf[\phi_1,\phi_2,\phi_3,\phi_4\rg]^\top$ and combining the equations~\eqref{eq:eq_set_1} and~\eqref{eq:eq_set_2} we obtain the system 
\begin{equation}
{\bf E}\boldsymbol\phi+\bol{\mathsf M}\bol\phi = \boldsymbol\phi^\inc \label{eq:system}
\end{equation}
where $\bol{\mathsf{ M}}$ is the $4\times 4$ block matrix integral  operator
$[\bol{\mathsf M}]_{i,j}:=\mathsf M_{i,j}$, $i,j=1,\ldots,4$,  
\begin{equation}\label{eq:data_def}
{\bf E}:=\begin{bmatrix}1\\ &\frac{1+\eta}{2}\\&&\gamma\\ &&&\gamma\end{bmatrix}\andtext \bol\phi^\inc :=\begin{bmatrix}f\\g\\0\\0\end{bmatrix}.
\end{equation}

Two observations about the system~\eqref{eq:system} are in order. The first one is that the last two equations in~\eqref{eq:system}, which account for the quasi-periodicity of the scattered field, need to be satisfied in all of $\R$. Being these equations as well as the associated density functions $\phi_3$ and $\phi_4$ defined in an unbounded interval, they need to be effectively truncated in order for them to be suitable to Nystr\"om or boundary element discretizations. We do so in the next section by means of the WGF  method. Secondly, note that the integral equation system~\eqref{eq:system} does not properly account for the radiation condition. Indeed, only the boundedness and the quasi-periodicity of the scattered field were used in its derivation. This important issue is also address in the next section. 

\begin{remark}
In light of~Proposition~\ref{prop:op_indeitities}, half of the operators~\eqref{eq:blocks2} can be significantly simplified in the case when $\Gamma_2$ and $\Gamma_3$ are parallel vertical lines. In fact, in such case we have
\begin{equation}\begin{aligned}
\mathsf M_{3,3}=&~-(1+\gamma^2)\mathsf K_1^{3,2},&\qquad \mathsf M_{3,4}=&~ (1-\gamma^2)\mathsf V^{3,2}_1,\smallskip\\\
\mathsf M_{4,3}=&~-(1-\gamma^2)\mathsf W^{2,3}_1,&\qquad  \mathsf M_{4,4}=&~(1+\gamma^2)\widetilde{\mathsf K}^{3,2}_1.
\end{aligned}\label{eq:blocks2_sim}\qed\end{equation}\end{remark}

\begin{remark} Note that other direct formulations can be used to account for the transmission conditions on~$\Gamma_1$. For instance, Kress-Roach equations~\eqref{eq:eq_set_1} can be replaced by the ones resulting from the well-known Costabel-Stephan formulation~\cite{costabel1985direct}, that can be easily derived by combining~\eqref{eq:gamma_1} and~\eqref{eq:trans_field} so as to eliminate $\phi_1$ and $\phi_2$ from the left-hand side of the equations. In this case we obtain 
\begin{equation*}
\sum_{q=1}^4\widetilde{\mathsf M}_{1,q}\phi_q=f \andtext\sum_{q=1}^4\widetilde{\mathsf M}_{2,q}\phi_q=g\quad\text{in}\quad [0,2\pi),
\label{eq:costabel}\end{equation*}
where
\begin{equation*}\begin{array}{llll}
\widetilde{\mathsf M}_{1,1}= -\mathsf K^{1,1}_2-\mathsf K^{1,1}_1,&\quad \widetilde{\mathsf M}_{1,2} =  \eta\mathsf V^{1,1}_1+\mathsf V^{1,1}_2,&\quad \widetilde{\mathsf M}_{1,3}=\mathsf M_{1,3},&\quad \widetilde{\mathsf M}_{1,4}=\mathsf M_{1,4},\smallskip\\
\widetilde{\mathsf M}_{2,1} =-\eta\mathsf W^{1,1}_2-\mathsf W^{1,1}_1,&\quad  \widetilde{\mathsf M}_{2,2}=\eta\mathsf{\widetilde K}^{1,1}_1+\eta\mathsf{\widetilde K}^{1,1}_2,&\quad\widetilde{\mathsf M}_{2,3}=\mathsf M_{2,3},&\quad\widetilde{\mathsf M}_{2,4}=\mathsf M_{2,4}.
\end{array}\label{eq:blocks_costable}\end{equation*} Unlike the advocated Kress-Roach approach, this formulation involves the (non-compact) hypersingular operator $\widetilde{\mathsf M}_{2,1}$ that negatively affect the conditioning of the discretized integral equation system, hindering the use of GMRES~\cite{saad1986gmres} and standard acceleration techniques based on fast matrix-vector products~\cite{rokhlin1990rapid}.~\qed
\end{remark}

\section{Windowed Green function method}
\label{sec:WGF}
In view of the definitions in~\eqref{eq:blocks2}, it is clear that several of the operators making up $\bol{\mathsf{M}}$ involve integration and evaluation over the unbounded curves $\Gamma_2$ or $\Gamma_3$. In order to reduce the BIE system~\eqref{eq:system} to a finite-size computational domain where standard BIE solvers can be applied, the domain of integration of the boundary integral operators over $\Gamma_2$ and $\Gamma_3$ has to be effectively truncated. We address this issue here by means of the WGF method~\cite{bruno2016windowed}.

The WGF method relies on the use of a  slow-rise window function $\chi(\cdot,cA,A)\in C^{\infty}_0(\R)$, $c\in(0,1)$, $A>0,$ which following~\cite{bruno2016windowed} is selected as
\begin{equation}
\chi(y,y_0,y_1) :=\left\{
\begin{array}{cl}
1&\text{if }|y|\leq y_0\\
\!\!\exp\left(\displaystyle\frac{2\e^{-1/u}}{u-1}\right)&\text{if } y_0<|y|<y_1, u=\frac{|y|-y_0}{y_1-y_0}\\
0&\text{if }|y|>y_1.
\end{array}\right.\label{eq:window_function}
\end{equation}

Note that  $\chi(\cdot,cA,A)$ vanishes together with all its derivatives in  $\{y\in\R:|y|> A\}$ and it equals one within  $\{y\in\R:|y|\leq   cA\}$.  In what follows we assume that $cA>\max\{h^+,-h^-\}$ so that the periodic array $D$ lies within the strip $\R\times[-cA,cA]$.

Next, letting
$$
w_{A}:=\chi(\cdot,cA,A)\circ\mathsf y_2\andtext w_A^c:=1-w_A
$$ 
 and replacing the split density 
$$\phi_j =  w_{A}\phi_j +w_{A}^c\phi_j\quad\text{for}\quad j=3,4$$ in~\eqref{eq:eq_set_1}-\eqref{eq:eq_set_2}, we obtain 
\begin{subequations}\begin{align}
\phi_1(t)+\sum_{q=1}^2 {\mathsf M}_{1,q}[\phi_q](t)+\sum_{q=3}^4 {\mathsf M}_{1,q}[w_A\phi_q](t)=&~f(t)-\psi_1(t),\ t\in[0,2\pi),\\
\lf(\frac{1+\eta}{2}\rg)\phi_2(t)+\sum_{q=1}^2 {\mathsf M}_{2,q}[\phi_q](t)+\sum_{q=3}^4 {\mathsf M}_{2,q}[w_A\phi_q](t)=&~g(t)-\psi_2(t),\ t\in[0,2\pi),\\
\gamma\phi_3(t)+\sum_{q=1}^2 {\mathsf M}_{3,q}[\phi_q](t)+\sum_{q=3}^4 {\mathsf M}_{3,q}[w_A\phi_q](t)=&~-\psi_3(t),\  t\in\R,\label{eq:win_param_1}\\
\gamma\phi_4(t)+\sum_{q=1}^2 {\mathsf M}_{4,q}[\phi_q](t)+\sum_{q=3}^4 {\mathsf M}_{4,q}[w_A\phi_q](t)=&~-\psi_4(t),\  t\in\R,\label{eq:win_param_2}
\end{align}\label{eq:win_eq_set_2}\end{subequations}
where the terms that were moved to the right-hand side  in~\eqref{eq:win_eq_set_2}  are the tail integrals 
\begin{equation}\label{eq:correction_functions}
\psi_p= {\mathsf M}_{p,3}[w_A^c\phi_3]+{\mathsf M}_{p,4}[w_A^c\phi_4],\quad p=1,\ldots,4.
\end{equation}

Our boundary integral equation formulation  relies on constructing suitable approximations of $\psi_p$, $p=1,\ldots,4,$ taking  into account the radiation condition~\eqref{eq:rad_cond_2} and the super-algebraic decay as $A\to\infty$ of certain oscillatory windowed integrals. Upon replacing $\psi_p$, $p=1,\ldots,4,$ by their respective approximations in~\eqref{eq:win_eq_set_2} and restricting  the integral equations~\eqref{eq:win_param_1} and ~\eqref{eq:win_param_2} to the bounded interval $[-A,A]$, we obtain a windowed integral equation suitable to be discretize by standard Nystr\"om or boundary element methods.

We then proceed to construct suitable approximations for the tail integrals $\psi_p$, $p=1,\ldots,4$. For the sake of presentation simplicity and without loss of generality in the remainder of this section we assume that ${\bf r}_2(t)=-\frac{L}2{\bf e}_1+t{\bf e}_2$ (i.e., $\mathsf y_2(t)=t$ ) for $t>|cA|$. From the general quasi-periodic expansion~\eqref{eq:general} of the scattered field it follows that within ${\rm supp}(w_A^c)=\{t\in\R:|t|\geq cA\}$ the parametrized traces $\phi_3$ and $\phi_4$~\eqref{eq:unknowns} associated with the unbounded curves $\Gamma_2$, can be expressed as
 \begin{equation}\begin{aligned}
\phi_3(t) =& \sum_{n\in\mathcal U\cup\mathcal V}\e^{-i\alpha_n\frac{L}2}\left\{B^\pm_n\e^{ \pm i\beta_n t }+C^\pm_n\e^{\mp i\beta_n t}\right\}+\sum_{n\in\mathcal W}\e^{-i\alpha_n\frac{L}{2}}\left\{B_n^\pm+C_n^\pm t\right\},\\
\phi_4(t) =& \sum_{n\in\mathcal U\cup\mathcal V}i\alpha_n\e^{-i\alpha_n\frac{L}2}\left\{B^\pm_n\e^{\pm i\beta_nt}+C^\pm_n\e^{\mp i\beta_n t}\right\}+\sum_{n\in\mathcal W}i\alpha_n\e^{-i\alpha_n\frac{L}{2}}\left\{B_n^\pm+C_n^\pm t\right\}
\end{aligned}\label{eq:Ray_traces}\end{equation}
for $\pm t>cA$. Splitting $w_A^c=1-w_A$ as $w_A^c= \chi_A^- +\chi_A^+$ where $\chi_A^- = {\bf 1}_{(-\infty,0)}w_A^c$ and $\chi_A^+={\bf 1}_{(0,\infty)}w_A^c$, and replacing~\eqref{eq:Ray_traces} in~\eqref{eq:correction_functions}, we arrive at
\begin{equation}\label{eq:approx_unbounded}
\psi_p =\psi_p^B + \psi_p^C,\quad p=1,\ldots,4, 
\end{equation}
where
\begin{equation}\begin{aligned}
\psi^B_p=  \sum_{n\in\mathbb Z}\e^{i\alpha_n\frac{L}2}B_n^+\lf\{{\mathsf M}_{p,3}\lf[\chi_A^+\e^{i\beta_n|\cdot|}\rg]+i\alpha_n{\mathsf M}_{p,4}\lf[\chi_A^+\e^{i\beta_n|\cdot|}\rg]\rg\}+\\
 \sum_{n\in\mathbb Z}\e^{i\alpha_n\frac{L}2}B_n^-\lf\{{\mathsf M}_{p,3}\lf[\chi_A^-\e^{i\beta_n|\cdot|}\rg]+i\alpha_n{\mathsf M}_{p,4}\lf[\chi_A^-\e^{i\beta_n|\cdot|}\rg]\rg\}
\end{aligned}\label{eq:1st_part}\end{equation}
and
\begin{equation}\begin{aligned}
\psi^C_p=&  \sum_{n\in\mathcal U\cup\mathcal V}\e^{-i\alpha_n\frac{L}2}C_n^+\lf\{{\mathsf M}_{p,3}\lf[\chi_A^+\e^{-i\beta_n|\cdot|}\rg]+i\alpha_n{\mathsf M}_{p,4}\lf[\chi_A^+\e^{-i\beta_n|\cdot|}\rg]\rg\}+\\
& \sum_{n\in\mathcal U\cup\mathcal V}\e^{-i\alpha_n\frac{L}2}C_n^-\lf\{{\mathsf M}_{p,3}\lf[\chi_A^-\e^{-i\beta_n|\cdot|}\rg]+i\alpha_n{\mathsf M}_{p,4}\lf[\chi_A^-\e^{-i\beta_n|\cdot|}\rg]\rg\}+\\
& \sum_{n\in\mathcal W}\e^{-i\alpha_n\frac{L}2}C_n^+\lf\{{\mathsf M}_{p,3}\lf[\chi_A^+\ \cdot\ \rg]+i\alpha_n{\mathsf M}_{p,4}\lf[\chi_A^+\ \cdot\ \rg]\rg\}+\\
& \sum_{n\in\mathcal W}\e^{-i\alpha_n\frac{L}2}C_n^-\lf\{{\mathsf M}_{p,3}\lf[\chi_A^-\ \cdot\ \rg]+i\alpha_n{\mathsf M}_{p,4}\lf[\chi_A^-\ \cdot\ \rg]\rg\}.
\end{aligned}\label{eq:2nd_part}\end{equation}

Let us first examine the term $\psi_p^B$, $p=1,\ldots,4$. In view of the boundedness of the Rayleigh coefficients  $B_n^\pm$ in~\eqref{eq:1st_part} and the exponential decay as $|t|\to\infty$ of the functions $\e^{i\beta_n|t|}$ for $\beta_n\in i\R_{>0}$ (i.e., $n\in\mathcal V$), we have that the approximation
\begin{equation}\begin{aligned}
\psi^B_p\approx\!\!\! \sum_{n\in\mathcal U\cup\mathcal W}\e^{i\alpha_n\frac{L}2}B_n^+\lf\{{\mathsf M}_{p,3}\lf[\chi_A^+\e^{i\beta_n|\cdot|}\rg]+i\alpha_n{\mathsf M}_{p,4}\lf[\chi_A^+\e^{i\beta_n|\cdot|}\rg]\rg\}+\\
\sum_{n\in\mathcal U\cup\mathcal W}\e^{i\alpha_n\frac{L}2}B_n^-\lf\{{\mathsf M}_{p,3}\lf[\chi_A^-\e^{i\beta_n|\cdot|}\rg]+i\alpha_n{\mathsf M}_{p,4}\lf[\chi_A^-\e^{i\beta_n|\cdot|}\rg]\rg\}
\end{aligned}\label{eq:1st_part_v1}\end{equation}
introduces errors that decrease exponentially fast as $A\to\infty$.

Next, for $\beta_n\in\R_{\geq0}$ (i.e., $n\in\mathcal U\cup\mathcal W$), we note that 
  ${\mathsf M}_{p,q}\lf[\chi_A^\pm \e^{ i\beta_n |\cdot|}\rg](t)$ for $p=1,2$, $t\in[0,2\pi)$, and for $p=3,4$, $t\in[-cA,cA]$, decays super-algebraically fast as $A\to\infty$ (i.e., faster than $O(((k_1+\beta_n)A)^{-m})$ for all $m\in\N$)~\cite{Bruno2015windowed,bruno2017windowed,perez2017windowed}. We refer the reader to Appendix~\ref{rem:sub_al_decay} for a detailed justification of this estimate. Given then the fast convergence of these windowed integrals as $A\to\infty$, we adopt to approximation 
  \begin{equation}\label{eq:approx_B}
  \psi_p^B\approx 0,\quad p=1,\dots,4.
\end{equation}

Let us now look into the terms $\psi_p^C$, $p=1,\ldots, 4$. In principle, the radiation condition~\eqref{eq:rad_cond_2} requires all the coefficients $C_n^\pm$ in~\eqref{eq:2nd_part} to vanish. These conditions can be easily incorporated in our formulation by simply setting $\psi_p^C=0$ which, together with~\eqref{eq:approx_B}, yield the following windowed integral equation:
\begin{equation}
{\bf{E}}\boldsymbol\phi_{\!A}+ \bol{\mathsf M}{\bf W}_{\!\!A}\boldsymbol\phi_{\!A} = \boldsymbol\phi^\inc\label{eq:window_system_2}
\end{equation}
where
\begin{equation}\label{eq:win_matrix}
{\bf W}_{\!\!A}(t) :=\begin{bmatrix}
1\\ &1\\ &&w_{\!A}(t)\\&&&w_{\!A}(t)\end{bmatrix},\quad t\in\R.
\end{equation}
Here, the first two equations of the system~\eqref{eq:window_system_2}, associated with the curve $\Gamma_1$, correspond to the parameter $t\in[0,2\pi)$, while the last two, associated with the truncated curve $\Gamma_{2,A}=\{\ner\in\R^2:\ner={\bf r}_2(t), |t|\leq A\}$, correspond to $t\in{\rm supp}(w_A) = [-A,A]$. Consequently, the entries $\phi_{j,A}$, $j=1,\dots,4$, of the solution vector $\bol \phi_A$ are considered functions $\phi_{j,A}:[0,2\pi]\to\C$ for $j=1,2$ and $\phi_{j,A}:[-A,A]\to\C$ for $j=3,4$. 

\section{An illustrative numerical example}\label{sec:illus_num_exam}
\begin{figure}[h!]
\centering	
\includegraphics[scale=1]{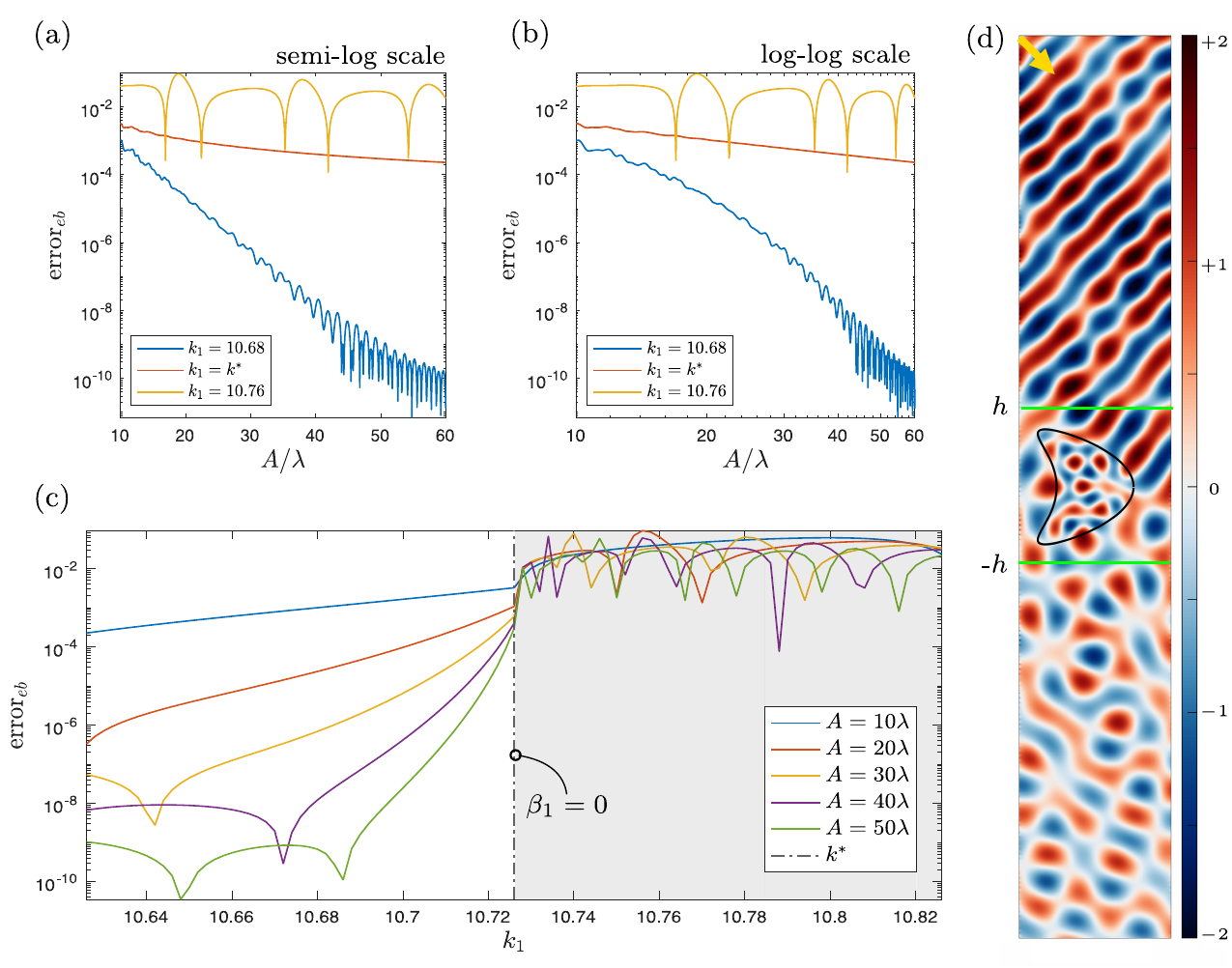}
\caption{Solution of the problem scattering of a planewave at $\theta^\inc=\frac{\pi}4$ by an infinite periodic array of kite-shaped obstacles obtained using the naive windowed BIE~\eqref{eq:window_system_2} and the scattered field approximation~\eqref{eq:approx_scat_fld_1} for $k_2=20$, $L=2$, $c=0.5$ and various window sizes $A$ and wavenumbers $k_1$ at and around a RW-anomaly configuration corresponding to $k_1=k^*\approx 10.7261$. Energy balance error~\eqref{eq:error_energy} as a function of~$A$ in (a)~semi-log  and (b)~log-log scale computed at $h=1$. (c) Wavenumber sweep of the energy balance error around~$k^*$. (d)  Real part of the computed total field within the region $[-\frac{L}2,\frac{L}2]\times[-cA,cA]$ for $k_1=10.68$ and $A=20\lambda$.}\label{fig:fail_example_1}
\end{figure}

In this section we consider a series of numerical experiments aimed at assessing the accuracy of the quasi-periodic problem solutions produced by the windowed integral equation~\eqref{eq:window_system_2}. In all such experiments we consider the diffraction and transmission of  a planewave~\eqref{eq:planewave} in TE polarization ($\eta=1$)  that impinges at an angle $\theta^\inc=\frac{\pi}4$ on an infinite array of period $L=2$ consisting of penetrable kite-shaped obstacles (see Figure~\ref{fig:fail_example_1}(d)) parametrized by 
\begin{equation}
{\bf r}_1(t)=\lf\{\frac12\cos t+\frac{13}{40} \cos 2 t-\frac{13}{40}\rg\}{\bf e}_1+ \frac34 \sin t\, {\bf e}_2, \quad t\in[0,2\pi).
\label{eq:kite_param}\end{equation}
For clarity of exposition the left ($\Gamma_2$) and right ($\Gamma_3$) hand side boundaries of the unit cell are selected  as straight vertical lines parametrized by
\begin{equation}\label{eq:param_gamma_2}
{\bf r}_2(t)=-\frac{L}2{\bf e}_1+t{\bf e}_2\andtext  {\bf r}_3(t)=\frac{L}2{\bf e}_1+t{\bf e}_2,\quad t\in \R,
\end{equation}
respectively. 

In our first experiment the error in the numerical solution is assessed by means of the energy balance relation (equation~\eqref{eq:conservation_energy} in Appendix~\ref{rem:energy}). We define the energy balance error as how much numerical solutions deviate from conserving energy, or more precisely, as
\begin{equation}\label{eq:error_energy}
\text{error}_{eb}:=\lf|2\real(\widetilde B_0^-)+\sum_{n\in\mathcal U}\frac{\beta_n}{\beta}\lf\{|\widetilde B_n^-|^2+|\widetilde B_n^+|^2\rg\}\rg|
\end{equation}
where the coefficients in~\eqref{eq:error_energy} are computed via (see Appendix~\ref{rem:energy})
\begin{equation}\label{eq:approx_coeff}
\widetilde B_n^\pm := \frac{\e^{\mp i\beta_n h}}{L}\int_{-\frac{L}2}^{\frac{L}2} u_A^s(x,\pm h)\e^{-i\alpha_n x}\de x
\end{equation}
using the WGF approximation of the scattered field given by
 \begin{equation}
u^s_A(\ner)  = (\mathcal D_1^1\phi_{A,1})(\ner)-\eta(\mathcal S_{1}^1\phi_{A,2})(\ner) +(\mathcal D_1^2-\gamma\mathcal D_1^3)[w_A\phi_{A,3}](\ner)-(\mathcal S_1^2-\gamma\mathcal S_1^3)[w_A\phi_{A,4}](\ner)\label{eq:approx_scat_fld_1}\end{equation}
for $\ner=(x,y)\in \Omega'$ with $\phi_{A,j}$, $j=1,\ldots,4$, denoting the components of the vector density $\bol\phi_A$  solution of~\eqref{eq:window_system_2}. Note that, as in the approximations that led to the windowed BIE system~\eqref{eq:window_system_2}, the errors produced by the windowed integrals in~\eqref{eq:approx_scat_fld_1} decay super-algebraically fast as $A\to\infty$ for $\ner\in\Omega_A' := \{(x,y)\in\Omega': \chi(y,cA,A)=1\}$ when the exact scattered field traces, $\phi_j$, $j=1,\ldots,4$ defined in~\eqref{eq:unknowns}, are used.

\begin{figure}[h!]
\centering	
 \includegraphics[scale=1]{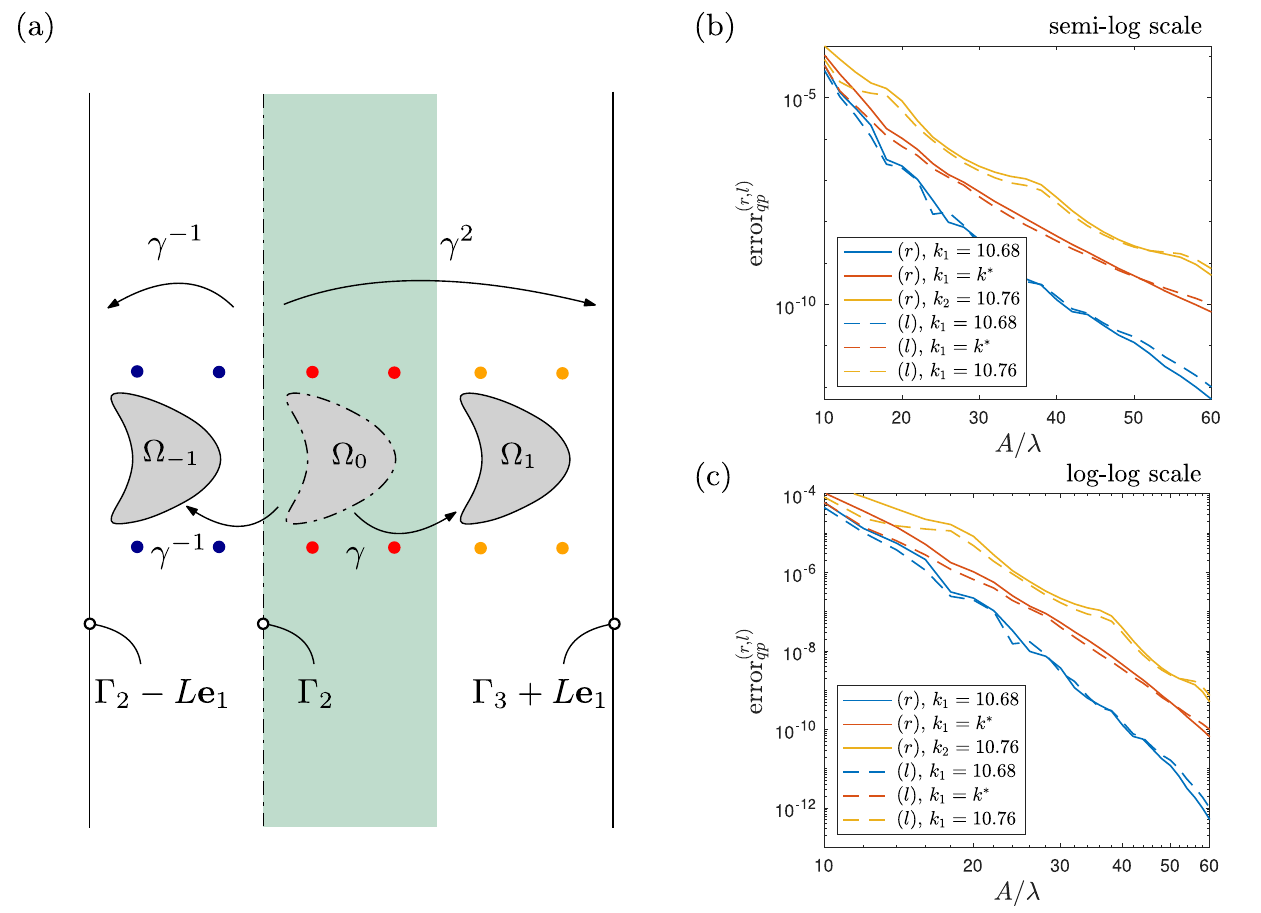}
\caption{Errors~\eqref{eq:qp_errors} in the quasi-periodicity condition of the numerical solution produced by the windowed integral equation~\eqref{eq:window_system_2}. (a) Depiction of the supercell configuration used to assess the left $(l)$ and right $(r)$ mismatch errors~\eqref{eq:qp_errors}. The density functions associated with the $3L$-periodic supercell are obtained from the densities of the middle $L$-periodic cell by multiplying them by $\gamma=\e^{i\alpha L}$ and $\gamma^{-1}=\e^{-i\alpha L}$ to transfer them from left to right and from right to left, respectively. Errors in semi-log (b) and log-log (c) scale for the exterior wavenumbers $k_1=10.68$, $k^*$, and $10.76$, and window sizes $A\in[10\lambda,60\lambda]$. }\label{fig:quasi_exp}
\end{figure}

Highly accurate numerical approximations of~$\bol\phi_A$ are used in all the examples presented in this section. These are obtained by numerically solving~\eqref{eq:window_system_2} by means of the spectrally accurate Martensen--Kussmaul (MK) Nystr\"om method~\cite[sec.~3.5]{COLTON:2012} employing a large number of discretization points (roughly, 8 points per wavelength on each of the relevant curves). The finite-domain integrals in~\eqref{eq:approx_coeff}, on the other hand, are computed using the trapezoidal quadrature rule which, by virtue of the fact that $u_A^s(\cdot,\pm h)\e^{-i\alpha_n \cdot}$ is smooth and approximately $L$-periodic (see~Figure~\ref{fig:quasi_exp}), it is expected to converge fast as the number of quadrature nodes increases. This choice of discretization methods and parameter values ensure that the dominant part of the energy balance error~\eqref{eq:error_energy} stems solely from the WGF approximation employed in~\eqref{eq:window_system_2} and~\eqref{eq:approx_scat_fld_1}. In what follows of this section we consider fixed parameter values $k_2=20$, $c=0.5$, and $h=1$.

Figures~\ref{fig:fail_example_1}(a)-(b) display the energy balance error~\eqref{eq:error_energy} in semi-log and log-log scale, respectively, as a function of the window size $A$ (measured in wavelengths $\lambda=2\pi/k_1$) for three different informative $k_1$ values at and around a RW-anomaly configuration corresponding to $k_1=k^*= 2\pi/(L(1-\sin\theta^\inc))\approx 10.7261$ ($\beta_1=0$ in this case). At $k_1=10.68$---to the left of $k^*$ when $\beta_1\in i\R_{>0}$---the expected super-algebraic convergence is achieved, as it can be seen in the blue nearly constant-slope  curve plotted in semi-log scale in Figure~\ref{fig:fail_example_1}(a). Figure~\ref{fig:fail_example_1}(d) displays the real part of the total field within the region $ [-\frac{L}2,\frac{L}2]\times[-cA,cA]=\overline{\Omega_A'\cup\Omega}$ for $A=20\lambda$, produced by the numerical evaluation of formulae~\eqref{eq:approx_scat_fld_1} and~\eqref{eq:GreenFormulaTrans}. At $k_1=k^*$ when $\beta_1=0$, in turn,  slow (algebraic) convergence is observed while at $k=10.76$---to the right of $k^*$ when $\beta_1\in \R_{>0}$---no convergence at all is observed. To examine this issue in more detail, a wider range of $k_1$ values is considered in Figure~\ref{fig:fail_example_1}(c), which shows a sweep of the error over the interval  $[k^*-0.1,k^*+0.1]$. Clearly, significant accuracy deterioration occurs at and around the RW-anomaly configuration for all the window sizes considered in this experiment.

There are two main factors that could explain the accuracy deterioration seen in Figure~\ref{fig:fail_example_1}(c). On the one hand we have the radiation condition, which is indirectly incorporated in our formulation by neglecting all the tail integrals $\psi^C_j$, $j=1,\ldots,4$,  in~\eqref{eq:win_eq_set_2}, and on the other hand, the quasi-periodicity condition~\eqref{eq:quasi-per} which is enforced through the equations~\eqref{eq:win_param_1} and~\eqref{eq:win_param_2} restricted to the interval $[-A,A]$.  

In order to verify the quasi-periodicity condition~\eqref{eq:quasi-per}, we consider the following experiment. First, the windowed integral equation~\eqref{eq:window_system_2} is solved using the MK method to obtain the approximate densities on $\Gamma_1$ and on the truncated vertical curves $\Gamma_{2,A}$ parametrized by ${\bf r}_2$ in~\eqref{eq:param_gamma_2} with $t$ restricted to $[-A,A]$. Then, assuming that the quasi-periodicity condition holds, we  ``transfer'' the densities to the boundaries of a $3L$-period supercell. Referencing to Figure~\ref{fig:quasi_exp}(a), we have that the supercell consists of the original obstacle's boundary $\Gamma_1=\p\Omega_0$, the shifted obstacles' boundaries $\Gamma_{1}-L{\bf e}_1:=\p\Omega_{-1}$ and $\Gamma_{1}+L{\bf e}_1 :=\p\Omega_1$, and the truncated parts $\Gamma_{2,A}-L{\bf e}_1$  and $\Gamma_{3,A}+L{\bf e}_1$ of the shifted vertical lines $\Gamma_{2}-L{\bf e}_1$  and $ \Gamma_{3}+L{\bf e}_1$ which are parametrized by ${\bf r}_2(\cdot)-L{\bf e}_1$ and  ${\bf r}_3(\cdot)+L{\bf e}_1$, respectively. Assuming that the quasi-periodicity condition holds, the densities associated with the supercell boundaries are: $\{\gamma^{-1}\phi_{A,1},\gamma^{-1}\phi_{A,2}\}$ on $\Gamma_{1}-L{\bf e}_1$, $\{\phi_{A,1},\phi_{A,2}\}$ on $\Gamma_1$, $\{\gamma\phi_{A,1},\gamma\phi_{A,2}\}$ on $\Gamma_{1}+L{\bf e}_1$,  $\{\gamma^{-1}\phi_{A,3},\gamma^{-1}\phi_{A,4}\}$ on $\Gamma_{2,A}-L{\bf e}_1$, and  $\{\gamma^2\phi_{A,3},\gamma^2\phi_{A,4}\}$ on $\Gamma_{3,A}+L{\bf e}_1$.  We then approximate the scattered field within the supercell as $u^s_A$ in~\eqref{eq:approx_scat_fld_1} but integrating on each of the relevant boundaries of the supercell using the aforementioned densities. To verify the quasi-periodicity condition we then introduce the right and left mismatch errors defined as
\begin{equation}\begin{aligned}
{\rm error}^{(r)}_{qp} :=&~\frac{\displaystyle\max_{p=1,\ldots,4} |u^s_A(\ner_p)-\gamma^{-1} u^s_A(\ner_p+L{\bf e}_1)|}{\displaystyle\max_{p=1,\ldots,4} |u^s_A(\ner_p)|}   \andtext\\
{\rm error}^{(l)}_{qp} :=&~\frac{\displaystyle\max_{p=1,\ldots,4} |u^s_A(\ner_p)-\gamma u^s_A(\ner_p-L{\bf e}_1)|}{\displaystyle\max_{p=1,\ldots,4} |u^s_A(\ner_p)|}
\end{aligned}\label{eq:qp_errors}\end{equation} respectively, 
where the sample points are ${\bf r}_1= (-0.5,-1)$, ${\bf r}_2= (0.5,-1)$, ${\bf r}_3= (-0.5,1)$, and ${\bf r}_4= (0.5,1)$ (they are depicted in Figure~\ref{fig:quasi_exp}(a) in red). The errors~\eqref{eq:qp_errors} corresponding to $k_1=10.68, k^*$ and $10.76$ are displayed in Figures~\ref{fig:quasi_exp}(b)-(c) in semi-log and log-log scales, respectively, for various window sizes $A\in [10\lambda,60\lambda]$. These results demonstrate that, although the enforcement of the quasi-periodicity condition deteriorates as $k_1$ approaches the RW-anomaly configuration, the mismatch errors still converge to zero super-algebraically fast as $A$ increases.
\begin{figure}[h!]
\centering	
\includegraphics[scale=0.48]{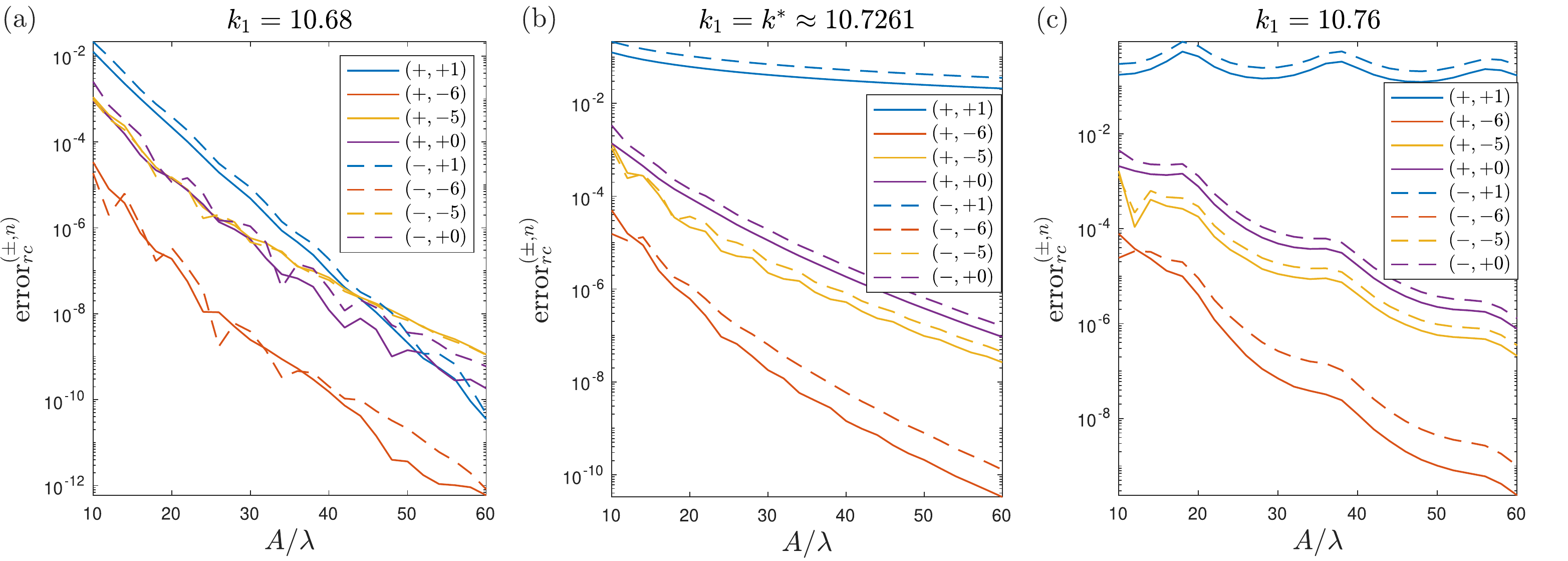}
\caption{Errors~\eqref{eq:error_rad} in the numerical solution obtained from the windowed integral equation~\eqref{eq:window_system_2} in the enforcement of the radiation condition~\eqref{eq:rad_cond_2}. Three different exterior wavenumbers are considered corresponding to $k_1=10.68$ in~(a), $k_1=k^*$ in~(b), and $k_1=10.76$ in~(c). The modes $n\in\mathcal C_{3k_1/4}$ used in these examples correspond to the smallest $\beta_n$ values arising in each case, which include $\beta_1$ that vanishes in the RW-tanomaly case $k_1=k^*$ in~(b).}\label{fig:fail_example_3}
\end{figure}

Since the quasi-periodicity condition does not seem to be the main factor that explains the poor convergence and the complete lack of it for certain wavenumbers~$k_1$, we are left with to examine the enforcement of the radiation condition. To do so, we introduce
\begin{equation}\label{eq:error_rad}\begin{split}
{\rm error}^{(+,n)}_{rc} :=&~ \lf|\frac{1}{L}\int_{-\frac{L}2}^{\frac{L}2}\left\{\p_y u^s_A(x, h)- i\beta_n u^s_A(x, h)\right\}\e^{-i\alpha_nx}\de x\rg|\andtext\\
{\rm error}^{(-,n)}_{rc} :=&~ \lf|\frac{1}{L}\int_{-\frac{L}2}^{\frac{L}2}\left\{\p_y u^s_A(x,- h)+ i\beta_n u^s_A(x,- h)\right\}\e^{-i\alpha_nx}\de x\rg|,
\end{split}\end{equation} where $u_A^s$ is defined in~\eqref{eq:approx_scat_fld_1}. These integrals are directly related to the amplitude $C_n^\pm$ of the non-radiative modes, which one expects to decrease super-algebraically fast as $A$ increases, thus measuring the correct enforcement of the radiation condition. 

Figure~\ref{fig:fail_example_3} displays the errors~\eqref{eq:error_rad} for $A\in[10\lambda,60\lambda]$, for the three representative wavenumbers $k_1=10.68, k^*$, and $10.76$, and for four modes~$n\in\mathcal C_{3k_1/4}=\{-6,-5,0,1\}$, where 
\begin{equation}\label{eq:reg_modes}
\mathcal C_\delta := \{n\in\Z: |\beta_n|\leq \delta\}.
\end{equation}
This set, which plays an important role below in Section~\ref{sec:corrected_BIE}, consists of the modes which are the closest to horizontally traveling waves. In the case $k_1=10.68$, which is considered in Figure~\ref{fig:fail_example_3}(a) and where $\beta_{-6}\approx 3.6844i$, $\beta_{-5}\approx6.8950$, $\beta_{0}\approx7.5519 $ and $\beta_{1}\approx0.5370i$, all the corresponding errors~\eqref{eq:error_rad} exhibit superalgebraic convergence as $A$ increases.  In turn, in the RW-anomaly case $k_1=k^*$, considered in Figure~\ref{fig:fail_example_3}(b) and where $\beta_{-6}\approx3.4429i$, $\beta_{-5}\approx7.0041$, $\beta_{0}\approx7.5845 $ and $\beta_1=0$,  slow convergence of ${\rm error}_{rc}^{(\pm,1)}$ is observed. Finally, in the case $k_1=10.76$, considered in Figure~\ref{fig:fail_example_3}(c) and where $\beta_{-6}\approx 3.2534i$, $\beta_{-5}\approx7.0835$, $\beta_{0}\approx7.6085 $ and $\beta_{1}\approx0.4624$,  we note that ${\rm error}_{rc}^{\pm,1}$ does not seem to converge at the all.  
These observations are consistent with the results displayed in  Figures~\ref{fig:fail_example_1}(a)-(b), that consider the overall energy balance error, and suggest that in practice the windowed BIE~\eqref{eq:window_system_2} on its own does not properly enforce the radiation condition of the problem. Indeed, the non-propagative modes corresponding to the smallest $\beta_n$ values, which are contained in $\mathcal C_\delta$, seem to be polluting the numerical solution. 


As it turns out, there is a subtle issue that explains the remarkable failure of the \emph{naive} windowed BIE~\eqref{eq:window_system_2}  for certain frequencies. In light of the estimates derived in Appendix~\ref{rem:sub_al_decay}, not only the tail integrals ${\mathsf M}_{p,q}\lf[\chi_A^\pm\e^{ i\beta_n|\cdot|}\rg]$  for $\beta_n\in\R_{> 0}$ decay super-algebraically fast as $A$ increases, but also  ${\mathsf M}_{p,q}\lf[\chi_A^\pm\e^{- i\beta_n|\cdot|}\rg]$  in~\eqref{eq:2nd_part}  as long as $\beta_n\in\R_{> 0}$ and $\beta_n\neq k_1$.  Indeed, these tend to zero faster than $O(((k_1-\beta_n)A)^{-m})$ for all $m\geq 1$ as $A\to\infty$. For a fixed $A>0$, this fact renders $C_n^\pm{\mathsf M}_{p,q}\lf[\chi_A^\pm\e^{- i\beta_n|\cdot|}\rg]$ for $\beta_n\in\R_{> 0}$, $\beta_n\neq k_1$, in~\eqref{eq:2nd_part} ``small'' regardless of the actual value of the coefficient $C_n^\pm$, hence in practice allowing the presence of non-radiative modes that eventually pollute the approximate solution of~\eqref{eq:window_system_2}. As it turns out, this is not much of an issue when ${\mathsf M}_{p,q}\lf[\chi_A^\pm\e^{- i\beta_n|\cdot|}\rg]$ converges slowly, i.e., when $\beta_n\approx k_1$,  but it certainly is  when ${\mathsf M}_{p,q}\lf[\chi_A^\pm\e^{- i\beta_n|\cdot|}\rg]$ converges fast, i.e., around a RW-anomaly configuration when $\beta_n\approx 0$. Indeed, this phenomenon explains why ${\rm error}_{rc}^{(\pm,1)}= \mp 2 i \beta_{1} \mathrm{e}^{i \beta_{1} h} C_1^\pm$  in Figure~\ref{fig:fail_example_3}(c), when $\beta_{1} \approx 0.4624$, does not seem to converge as $A$ increases, while in turn ${\rm error}_{rc}^{(\pm,0)}$ when $\beta_{0} \approx 7.6085$, exhibits fast convergence. Interestingly, this phenomenon is present even in connection with the divergent tail integrals in~\eqref{eq:2nd_part} corresponding to ${\mathsf M}_{p,q}\lf[\chi_A^\pm\,\cdot\,\rg]$ for $\beta_n=0$ and  ${\mathsf M}_{p,q}\lf[\chi_A^\pm\e^{- i\beta_n|\cdot|}\rg]$ for $\beta_n\in i\R_{ >0}$ and $\beta_n\approx 0$, due to the slow divergence of the complementary integrals along the bounded interval $[-A,A]$. This is for instance observed in Figure~\ref{fig:fail_example_3}(b) which shows the slow convergence of ${\rm error}_{rc}^{(\pm,1)}=  C_1^\pm$.

\section{Corrected windowed integral equation}\label{sec:corrected_BIE}

This section presents a correction to the windowed integral equation~\eqref{eq:window_system_2} that renders it robust for all frequencies. To do so we first examine the case $\mathcal W=\emptyset$, in which $\beta_n\neq0$ for all $n\in\Z$, while the RW-anomaly case $\mathcal W\neq \emptyset$ is dealt with in Section~\ref{sec:wood_anom} below.

In order to address the issue discussed in the previous section, we proceed to retain certain critical coefficients $C_n^\pm$ as unknown variables in our formulation. Guided by the results of the previous section, such coefficients are selected as $C_n^\pm$ for $n\in \mathcal C_\delta$ defined in~\eqref{eq:reg_modes} where $\delta>0$ is a user-provided parameter.  The conditions $C_m^\pm=0$ for $n\in\mathcal C_\delta$ will be then indirectly enforced through the integral form of the radiation condition~\eqref{eq:rad_cond_2}.

Doing so the tail integrals~\eqref{eq:2nd_part} will no longer be set to zero but become
\begin{equation}\begin{aligned}
\psi^C_p\approx&  \sum_{n\in\mathcal C_\delta}\{C_n^+\Psi_{n,p}^++C_n^-\Psi_{n,p}^-\},\quad p=1,\ldots,4,
\end{aligned}\label{eq:3rd_part}\end{equation}
where the functions $\Psi^\pm_{n,p}$  are (formally) defined as
\begin{equation}
\Psi^\pm_{n,p} =
\e^{-i\alpha_n\frac{L}2}\lf\{{\mathsf M}_{p,3}\lf[\chi_A^\pm\e^{\mp i\beta_n|\cdot|}\rg]+i\alpha_n{\mathsf M}_{p,4}\lf[\chi_A^\pm\e^{\mp i\beta_n|\cdot|}\rg]\rg\},\quad n\in \mathcal C_\delta.
\label{eq:crr_terms}\end{equation}

Note that~\eqref{eq:crr_terms} involves evaluation of improper integrals over the unbounded intervals $(-\infty,-cA]$ and $[cA,\infty)$ that either cannot be evaluated in closed form or simply diverge. Therefore,  computable approximations of $\Psi_{n,p}^\pm$  in~\eqref{eq:crr_terms} are needed. Here, we resort to Green's representation formula~\eqref{eq:GreenFormula2} and  the fact that
\begin{equation}\label{eq:neg_approx}
{\mathsf M}_{p,3}\lf[\chi_A^\mp\e^{\mp i\beta_n|\cdot|}\rg]+i\alpha_n{\mathsf M}_{p,4}\lf[\chi_A^\mp\e^{\mp i\beta_n|\cdot|}\rg],\quad n\in\Z,
\end{equation}
tend to zero either super-algebraicaly (for $n\in\mathcal U$) or exponentially  (for $n\in\mathcal V$) fast as $A$ increases (see Appendix~\ref{rem:sub_al_decay}), to produce such approximations.

Let us first consider the case  $n\in\mathcal U\cup\{m\in\Z:\beta_m\neq k_1\}$ ($\beta_n\in \R_{>0}$, $\beta_n\neq k_1$) for which~\eqref{eq:crr_terms} are well-defined conditionally convergent integrals. Approximations of $\Psi^\pm_{n,p}$ for the remaining $\beta_n$ values are obtained by simply considering the analytical extension of the resulting expressions that depend smoothly on $\beta_n$.

Using then the fact that~\eqref{eq:neg_approx} becomes negligible for large $A$ values, we can add it to $\Psi^\pm_{n,p}$ in~\eqref{eq:crr_terms} to form 
\begin{align*}
\Psi^\pm_{n,p} \approx&~ \e^{-i\alpha_n\frac{L}2}\lf\{{\mathsf M}_{p,3}\lf[w_A^c\e^{\mp i\beta_n|\cdot|}\rg]+i\alpha_n{\mathsf M}_{p,4}\lf[w_A^c\e^{\mp i\beta_n|\cdot|}\rg]\rg\},
\end{align*}
where we used the identities $w_A^c=1-w_A=\chi_A^\pm+\chi_A^\mp$. Then, introducing the notation
\begin{equation}\begin{aligned}
\phi_{n,1}^{\pm} := u_n^{\pm}\big|_{\Gamma_1}\circ {\bf r}_1=\e^{i\alpha_n\mathsf x_1\pm i\beta_n\mathsf y_1},\qquad  \phi_{n,2}^{\pm} :=\p_n u_n^\pm\big|_{\Gamma_1}\circ {\bf r}_1=i{\bf n}_1\cdot(\alpha_n,\pm\beta_n)\e^{i\alpha_n\mathsf x_1\pm i\beta_n\mathsf y_1},\\
\phi_{n,3}^{\pm} := u_n^{\pm}\big|_{\Gamma_2}\circ {\bf r}_2=\e^{i\alpha_n\mathsf x_2\pm i\beta_n\mathsf y_2},\qquad \phi_{n,4}^{\pm} :=\p_n u_n^\pm\big|_{\Gamma_2}\circ {\bf r}_2=i{\bf n}_2\cdot(\alpha_n,\pm\beta_n)\e^{i\alpha_n\mathsf x_2\pm i\beta_n\mathsf y_2},
\end{aligned}\label{eq:traces_modes}\end{equation}
for the parametrized traces of the Rayleigh modes~\eqref{eq:prog_modes}, and exploiting the linearity of the integral operators~${\mathsf M}_{p,q}$, we arrive at
\begin{equation}\begin{aligned}
\Psi^\pm_{n,p} \approx&~{\mathsf M}_{p,3}\lf[(1-w_A)\phi_{n,3}^\mp\rg]+ {\mathsf M}_{p,4}\lf[(1-w_A)\phi_{n,4}^\mp\rg]=\\
&-\Phi_{n,p}^\pm-\lf\{{\mathsf M}_{p,3}\lf[w_A\phi_{n,3}^\mp\rg]+ {\mathsf M}_{p,4}\lf[w_A\phi_{n,4}^\mp\rg]\rg\}
\end{aligned}\label{eq:comp_approx_win}\end{equation}
where closed-form expressions for the functions
\begin{equation*}
\Phi_{n,p}^\pm = -{\mathsf M}_{p,3}\phi_{n,3}^{\mp}-{\mathsf M}_{p,4}\phi_{n,4}^\mp
\end{equation*} can be found.
Indeed, from the definition of the operators ${\mathsf M}_{p,q}$, $p=1,\ldots,4$ and $q=2,3,$ in~\eqref{eq:blocks1} and~\eqref{eq:blocks2}, and Green's representation formula~\eqref{eq:GreenFormula2}, we find that 
\begin{equation}\begin{split}
\Phi_{n,p}^\pm =&-\begin{cases}(\gamma \mathsf{K}_{1}^{1,3}-\mathsf{K}_{1}^{1,2}) \phi_{n,3}^\mp-(\gamma \mathsf{V}_{1}^{1,3}-\mathsf{V}_{1}^{1,2})\phi_{n,4}^\mp&p=1\\
(\gamma \mathsf{W}_{1}^{1,3}-\mathsf{W}_{1}^{1,2}) \phi_{n,3}^\mp-(\gamma \widetilde{\mathsf{K}}_{1}^{1,3}-\widetilde{\mathsf{K}}_{1}^{1,2})\phi_{n,4}^\mp&p=2\\
(\gamma^{2} \mathsf{K}_{1}^{2,3}-\mathsf{K}_{1}^{3,2}) \phi_{n,3}^\mp-(\gamma^{2}\mathsf V_{1}^{2,3}-\mathsf V_{1}^{3,2})\phi_{n,4}^\mp&p=3\\
(\gamma^{2} \mathsf{W}_{1}^{2,3}-\mathsf{W}_{1}^{3,2}) \phi_{n,3}^\mp-(\gamma^{2} \widetilde{\mathsf K}_{1}^{2,3}-\widetilde{\mathsf K}_{1}^{3,2})\phi_{n,4}^\mp&p=4
\end{cases}\\
=& \begin{cases} \phi_{n,p}^\mp&p=1,2\\
\gamma \phi_{n,p}^\mp&p=3,4.
\end{cases}\end{split}\label{eq:corrections}
\end{equation}

For $n\in\mathcal U\cup\{m\in\Z:\beta_m\neq k_1\}$ we have hence produced a computable approximation~\eqref{eq:comp_approx_win} of the modal integrals $\Psi^\pm_{n,p}$~\eqref{eq:crr_terms} with errors that decay super-algebraically fast as the window size $A$ increases.  Such an approximation consists of the closed-form expression~\eqref{eq:corrections} and the finite-domain windowed integrals in~\eqref{eq:comp_approx_win} that can be evaluated numerically. Corresponding computable expressions for $\Psi^\pm_{n,p}$ in the case  $n\in\mathcal V\cup\{m\in\Z:\beta_m=k_1\}$ are obtained by analytically extending the formula on the right-hand side of~\eqref{eq:comp_approx_win} to $\beta_n$ values.

We are now in position to write the corrected windowed BIE in the case $\mathcal W=\emptyset$. Letting 
\begin{equation}
\bol \Psi_n^\pm =-\begin{bmatrix}\phi_{n,1}^\mp\\\phi_{n,2}^\mp\\\gamma\phi_{n,3}^\mp\\\gamma\phi_{n,4}^\mp\end{bmatrix}-\bol{\mathsf M}{\bf W}_{\!\!A}\begin{bmatrix}0\\0\\ \phi^\mp_{n,3}\\ \phi^\mp_{n,4}\end{bmatrix}
\label{eq:vec_aux_func}\end{equation} and using~\eqref{eq:comp_approx_win} we obtain that the BIE can be expressed in vector form as
\begin{equation}
{\bf E}\boldsymbol\phi_{\!A}+ \bol{\mathsf M}{\bf W}_{\!\!A}\boldsymbol\phi_{\!A} +\sum_{n\in\mathcal C_\delta}\lf\{C_n^+\bol\Psi_n^++C_n^-\bol\Psi_n^-\rg\}= \bol\phi^\inc\label{eq:window_system_3}
\end{equation}
where, as~\eqref{eq:window_system_2}, the first two equations hold in the interval $[0,2\pi)$ while the last two hold in $[-A,A]$. 

The system~\eqref{eq:window_system_3} together with the radiation condition~\eqref{eq:rad_cond_2}, which yields the equations
\begin{equation}\label{eq:coeff_eqns}
C_n^\pm = \mp\frac{\e^{ i\beta_n h}}{2i\beta_n L}\int_{-\frac{L}2}^{\frac{L}2}\left\{\p_y u^s_A(x,\pm h)\mp i\beta_n u^s_A(x,\pm h)\right\}\e^{-i\alpha_nx}\de x=0,\quad  n\in \mathcal C_\delta,
\end{equation}
that couple the coefficients $C_n^\pm$,  $n\in\mathcal C_\delta$, with the unknown density vector $\bol\phi_{\!A}$ through a suitable WGF approximation $u^s_A$ of the scattered field,  form the system that we solve numerically. The parameter $h>0$ in~\eqref{eq:coeff_eqns} has to satisfy $\max\{h^+,-h^-\}<h<cA$.

The next section presents a suitable WGF approximation $u_A^s$ of the scattered field to be used in~\eqref{eq:coeff_eqns}, which unlike the naive approximation~\eqref{eq:approx_scat_fld_1}, explicitly accounts for the non-radiative modes in~\eqref{eq:3rd_part}.

\subsection{Windowed Green function approximation of the scattered field}
 As the matrix integral operator in~\eqref{eq:system}, the representation formula~\eqref{eq:GF_scat} of the scattered field involves the computation of layer potentials along the unbounded curve $\Gamma_2$. In view of the discussion of the previous section we proceed to utilize the following approximation of the $\Gamma_2$ traces of $u^s$:
\begin{equation}\label{eq:approx_gamma_2}
\phi_j \approx w_A\phi_{\!A,j}  + \sum_{n\in\mathcal C_\delta} \lf\{C_n^+\chi_A^+\phi_{n,j}^-+C_n^-\chi_A^-\phi_{n,j}^+\rg\},\quad j=3,4,
\end{equation}
where $\phi_{\!A,j}$, $j=1,\ldots,4,$ denote the entries of the vector density $\bol\phi_{\!A}$ in~\eqref{eq:window_system_2} and $\phi_{n,j}^+$, $j=1,\ldots,4$ denote the traces of the (non-radiative) modes introduced in~\eqref{eq:traces_modes}.  The presence of such modes in~\eqref{eq:approx_gamma_2} accounts for the fact that the integral equation~\eqref{eq:window_system_2} as well as the integral representation of the scattered field~\eqref{eq:approx_scat_fld_1} used in Section~\ref{sec:illus_num_exam} do not properly account for the radiation condition. 

 Replacing~\eqref{eq:approx_gamma_2} in the integral representation of the scattered field~\eqref{eq:GF_scat} we obtain
\begin{equation}\begin{split}
u^s(\ner)  \approx&~ (\mathcal D_1^1\phi_{A,1})(\ner)-\eta(\mathcal S_{1}^1\phi_{A,2})(\ner) +\\
&~ (\mathcal D_1^2-\gamma\mathcal D_1^3)[w_A\phi_{A,3}](\ner)-(\mathcal S_1^2-\gamma\mathcal S_1^3)[w_A\phi_{A,4}](\ner)+\\
&\sum_{n\in\mathcal C_\delta}C_n^+\lf\{ (\mathcal D_1^2-\gamma\mathcal D_1^3)[\chi_A^+\phi^-_{n,3}](\ner)-(\mathcal S_1^2-\gamma\mathcal S_1^3)[\chi_A^+\phi^-_{n,4}](\ner)\rg\}+\\
&\sum_{n\in\mathcal C_\delta}C_n^-\lf\{ (\mathcal D_1^2-\gamma\mathcal D_1^3)[\chi_A^-\phi^+_{n,3}](\ner)-(\mathcal S_1^2-\gamma\mathcal S_1^3)[\chi_A^-\phi^+_{n,4}](\ner)\rg\},\quad \ner\in \Omega'\end{split}\label{eq:approx_scat_fld}\end{equation}
where the layer potentials are defined in~\eqref{eq:param_pots}.

To produce a computable approximation of the modal terms in~\eqref{eq:approx_scat_fld} we resort to the above mentioned properties of windowed oscillatory integrals to note that, for a target point $\ner\in\Omega'_A = \{\ner=(x,y)\in\Omega: \chi(y,cA,A)=1\}$, the integrals 
\begin{equation}\label{eq:approx_pot}
(\mathcal D_1^2-\gamma\mathcal D_1^3)[\chi_A^\pm\phi^\mp_{n,3}](\ner)-(\mathcal S_1^2-\gamma\mathcal S_1^3)[\chi_A^\pm\phi^\mp_{n,4}](\ner)
\end{equation}
can be effectively approximated by 
$$
(\mathcal D_1^2-\gamma\mathcal D_1^3)[w^c_A\phi^\mp_{n,3}](\ner)-(\mathcal S_1^2-\gamma\mathcal S_1^3)[w_A^c\phi^\mp_{n,4}](\ner) 
$$
with errors 
$$
(\mathcal D_1^2-\gamma\mathcal D_1^3)[\chi_A^\mp\phi^\mp_{n,3}](\ner)-(\mathcal S_1^2-\gamma\mathcal S_1^3)[\chi_A^\mp\phi^\mp_{n,4}](\ner) 
$$ that converge to zero either super-algebraically fast for  $n\in\mathcal U$ (i.e., $\beta_n\in \R_{>0}$) or exponentially fast for $n\in\mathcal V$ (i.e., $\beta_n\in i\R_{>0}$) as $A\to\infty$. 

Therefore, letting  $\phi_{A,j}^c$, $j=1,\ldots,4$, denote the entries of the corrected vector density 
\begin{equation}\label{eq:corrected_density}
\bol\phi_A^c=\bol\phi_A-\sum_{n\in\mathcal C_\delta}\lf\{C_n^+ \begin{bmatrix}0\\0\\ \phi^-_{n,3}\\ \phi^-_{n,4}\end{bmatrix}+C_n^- \begin{bmatrix}0\\0\\ \phi^+_{n,3}\\ \phi^+_{n,4}\end{bmatrix}\rg\}
\end{equation}
we define our WGF approximation of the scattered field as
\begin{equation}\label{eq:approx_scat_fld_2}\begin{split}
u^s_A(\ner)  =&~ (\mathcal D_1^1\phi^c_{A,1})(\ner)-\eta(\mathcal S_{1}^1\phi^c_{A,2})(\ner) + (\mathcal D_1^2-\gamma\mathcal D_1^3)[w_A\phi^c_{A,3}](\ner)-\\
&(\mathcal S_1^2-\gamma\mathcal S_1^3)[w_A\phi^c_{A,4}](\ner)+\sum_{n\in\mathcal C_\delta}\lf\{C_n^+u_n^-(\ner)+C_n^{-}u_n^+(\ner)\rg\}
\end{split}\end{equation}
for $\ner\in \Omega'_A$ where  the last two terms were obtained by direct application of Green's representation formula~\eqref{eq:GreenFormula2}. 

With this expression at hand, we can now easily incorporate the conditions~\eqref{eq:coeff_eqns} into the integral equation system. In order to do so we define the functionals:
\begin{equation}\label{eq:L_functionals}\begin{split}
\mathsf L_n^\pm\bol\phi=\frac{1}{L} \int_{-\frac{L}2}^{\frac{L}2}&\lf[\mathsf (\p_y\mathcal D_1^1\mp i\beta_n\mathcal D_1^1)\phi_1-\eta(\p_y\mathcal S_{1}^1\mp i\beta_n\mathcal S_{1}^1)\phi_2\rg.+
\lf. \{\p_y\mathcal D_1^2\mp i\beta_n\mathcal D_1^2 -\gamma(\p_y \mathcal D_1^3\mp i\beta_n\mathcal D_1^3)\}\phi_{3}\rg.-\\
&~\lf.\{\p_y\mathcal S_1^2\mp i\beta_n\mathcal S_1^2-\gamma(\p_y\mathcal S_1^3\mp i\beta_n \mathcal S_1^3)\}\phi_{4}\rg]({\bf r}_h^\pm(t))\e^{-i\alpha_nt}\de t
\end{split}\end{equation}
where ${\bf r}_h^\pm(t) =\pm h{\bf e}_2+t{\bf e}_1$, with which conditions~\eqref{eq:coeff_eqns}  using~\eqref{eq:approx_scat_fld_2} can be readily expressed as
\begin{equation}\begin{aligned}
 C_n^+=\frac{\e^{i\beta_n h}}{2i\beta_n}\mathsf L_{n}^+\lf[{\bf W}_{\!\!A}\bol\phi_A^c\rg]\andtext
  C_n^-=- \frac{\e^{i\beta_n h}}{2i\beta_n}\mathsf L_{n}^-\lf[{\bf W}_{\!\!A}\bol\phi_A^c\rg],\quad n\in\mathcal C_{\delta} 
\end{aligned}\end{equation}
(note that we are still assuming that $\beta_n\neq 0$ for all $n\in\mathcal C_\delta$, i.e., $\mathcal W=\emptyset$).

Therefore, both~\eqref{eq:window_system_3} and~\eqref{eq:coeff_eqns} can be recast as a single corrected windowed BIE system:
\begin{equation}
\bf {E}\boldsymbol\phi^c_{\!A}+ \bol{\mathsf M}_c{\bf W}_{\!\!A}\boldsymbol\phi^c_{\!A} = \boldsymbol\phi^\inc\label{eq:corrected_lin_sys}
\end{equation}
 for the corrected vector density~\eqref{eq:corrected_density}, where letting $ \bol\Phi_n^\pm=\begin{bmatrix}\phi_{n,1}^\mp\\\phi_{n,2}^\mp\\0\\0\end{bmatrix} 
$, the corrected matrix operator is given by
\begin{equation}\label{eq:corrected_mat_v1}
\bol{\mathsf M}_c=\bol{\mathsf M}+\sum_{n\in\mathcal C_\delta}\frac{\e^{i\beta_n h}}{2i\beta_n}\lf\{\bol\Phi_n^-\mathsf  L_{n}^--\bol\Phi_n^+\mathsf L_{n}^+\rg\}
\end{equation}
in the case $\mathcal W=\emptyset$.
\begin{remark} Note that the functionals~$\mathsf L_n^\pm$ defined in \eqref{eq:L_functionals} entail evaluation of singular integrals. This is so because the layer potentials $\mathcal D^i_1$ and $\mathcal S_1^i$ involve integration along the unit-cell boundaries $\Gamma_i$, $i=2,3$, which are intersected by the horizontal line segments parametrized by ${\bf r}^\pm_h$. 

To avoid  this issue altogether we leverage the quasi-periodicity condition satisfied by the scattered field and express it by means of Green's representation formula applied within a three-period wide cell, like the one employed in the numerical examples of Figure~\ref{fig:quasi_exp}. The scattered field is then produced through integration on the super-cell walls  $\Gamma_2-L{\bf e}_1$ and $\Gamma_3+L{\bf e}_1$, which are parametrized by ${\bf r}_2(\cdot)-L{\bf e}_1$ and ${\bf r}_2(\cdot)+2L{\bf e}_1$, respectively, as well as on the annexed left and right obstacle boundaries $\Gamma_{1}-L{\bf e}_1$ and $\Gamma_{1}+L{\bf e}_1$, which are  parametrized by ${\bf r}_1(\cdot)-L{\bf e}_1$ and ${\bf r}_1(\cdot)+L{\bf e}_1$, respectively.  The densities on the new curves are given by multiplying the original densities by $\gamma^{-1}$ and $\gamma$ depending on whether the new curve corresponds to left or right $L$-translation of the original curve, respectively. Doing so the functionals can be recast as
\begin{equation}\label{eq:L_functionals_v2}\begin{split}
\mathsf L_n^\pm\bol\phi=&\frac{1}{L} \int_{-\frac{L}2}^{\frac{L}2}\lf[\lf\{\p_y(\mathcal D_1^1+\gamma^{-1}\mathcal D_1^{1-L}+\gamma \mathcal D_1^{1+L})\mp i\beta_n(\mathcal D_1^1+\gamma^{-1}\mathcal D_1^{1-L}+\gamma \mathcal D_1^{1+L})\rg\}\phi_1\rg.-\\
&~\eta\lf\{(\p_y(\mathcal S_{1}^1+\gamma^{-1}\mathcal S_{1}^{1-L}+\gamma \mathcal S_{1}^{1+L})\mp i\beta_n(\mathcal S_{1}^1+\gamma^{-1}\mathcal S_{1}^{1-L}+\gamma \mathcal S_{1}^{1+L})\rg\}\phi_2+\\
&~\lf. \lf\{\gamma^{-1}(\p_y{\mathcal D}_1^{2-L}\mp i\beta_n{\mathcal D}_1^{2-L}) -\gamma^2(\p_y {\mathcal D}_1^{3+L}\mp i\beta_n{\mathcal D}_1^{3+L})\rg\}\phi_{3}\rg.-\\
&~\lf.\lf\{\gamma^{-1}(\p_y{\mathcal S}_1^{2-L}\mp i\beta_n{\mathcal S}_1^{2-L})-\gamma^2(\p_y\mathcal S_1^{3+L}\mp i\beta_n \mathcal S_1^{3+L})\rg\}\phi_{4}\rg]({\bf r}_h^\pm(t))\e^{-i\alpha_nt}\de t
\end{split}\end{equation}
in terms of the layer potentials: $\mathcal D_1^{1+L}$ and  $\mathcal S_1^{1+L}$ associated with $\Gamma_{1}+L{\bf e}_1$; $\mathcal D^{1-L}_1$ and $\mathcal S^{1-L}_1$ associated with $\Gamma_{1-L}$;  ${\mathcal D}_{1}^{2-L}$ and ${\mathcal S}_1^{2-L}$ associated with $\Gamma_2-L{\bf e}_1$; and, ${\mathcal D}_1^{3+L}$ and ${\mathcal S}_1^{3+L}$ associated with $\Gamma_3+L{\bf e}_1$.~\qed
\end{remark}

\subsection{Corrected windowed integral equation at Rayleigh-Wood anomalies}\label{sec:wood_anom}

Finally, in order to extend~\eqref{eq:corrected_lin_sys} to the challenging RW-anomaly case, i.e., when $\beta_n=0$ for some $n\in\mathcal C_\delta$ ($\mathcal W\neq\emptyset$), we resort to  L'H\^opital's rule. In detail, we evaluate the correcting terms in~\eqref{eq:corrected_mat_v1} associated with $n\in \mathcal W$ as the limit
\begin{equation}\begin{aligned}
\lim_{\beta_n\to 0}\frac{\e^{i\beta_n h}}{\beta_n}\lf\{\bol\Phi_n^-\mathsf  L_{n}^--\bol\Phi_n^+\mathsf L_{n}^+\rg\}=&~
\p_{\beta_n}\lf\{\bol\Phi_n^-\mathsf  L_{n}^--\bol\Phi_n^+\mathsf L_{n}^+\rg\}\Big|_{\beta_n=0}\\
=&~\lf\{\p_{\beta_n}\bol\Phi_n^-\big|_{\beta_n=0}\mathsf  L_{n}^-\big|_{\beta_n=0}+\bol\Phi_n^-\big|_{\beta_n=0}\p_{\beta_n}\mathsf  L_{n}^-\big|_{\beta_n=0}\rg\}-\\
&~\lf\{\p_{\beta_n}\bol\Phi_n^+\big|_{\beta_n=0}\mathsf  L_{n}^+\big|_{\beta_n=0}+\bol\Phi_n^+\big|_{\beta_n=0}\p_{\beta_n}\mathsf  L_{n}^+\big|_{\beta_n=0}\rg\}.
\end{aligned}\end{equation}
Doing so the general expression for the corrected matrix operator in~\eqref{eq:corrected_mat_v1} becomes 
\begin{equation}\label{eq:corrected_mat}
\bol{\mathsf M}_c:=\bol{\mathsf M}+\frac{1}{2i}\sum_{n\in\mathcal C_\delta\setminus\mathcal W}\frac{\e^{i\beta_n h}}{\beta_n}\lf\{\bol\Phi_n^-\mathsf  L_{n}^--\bol\Phi_n^+\mathsf L_{n}^+\rg\}+\frac{1}{2i}\sum_{n\in\mathcal W}\lf\{\bol\Phi_n\p_{\beta_n}(\mathsf  L_{n}^--\mathsf L_{n}^+)+\p_{\beta_n}\bol\Phi_n(\mathsf  L_{n}^-+\mathsf L_{n}^+)\rg\}
\end{equation}
where we have introduced the vectors $\bol\Phi_n:=\bol\Phi_n^\pm$ and  $$\p_{\beta_n}\bol\Phi_n=\begin{bmatrix} i\mathsf y_1\\{\bf n}_1\cdot(-\mathsf y_1\alpha_n,i)\\0\\0\end{bmatrix}\e^{i\alpha_n x}$$ for $n\in\mathcal W$, which  correspond to the $\Gamma_1$ traces of the Raleigh modes $u_n$ and $iv_n$ defined in~\eqref{eq:modes}. The $\beta_n$-derivative of the functionals $\mathsf L_n^\pm$ are given by 
\begin{equation}\label{eq:L_functionals__der}\begin{split}
\p_{\beta_n}\mathsf L_n^\pm\bol\phi=\mp\frac{i}{L} &\int_{-\frac{L}2}^{\frac{L}2}\lf[  \mathcal D_1^1\phi_1-\eta \mathcal S_{1}^1\phi_2+ \gamma\mathcal D_1^{1+L}\phi_1-\gamma\eta \mathcal S_{1}^{1+L}\phi_2+\rg.\\
&\gamma^{-1} \mathcal D_1^{1-L}\phi_1-\gamma^{-1}\eta \mathcal S_{1}^{1-L}\phi_2+ \gamma^{-1} {\mathcal D}_1^{2-L}\phi_{3} -\gamma^2 {\mathcal D}_1^{3+L}\phi_{3}-\\
&~\lf.\gamma^{-1} {\mathcal S}_1^{2-L}\phi_{4}+\gamma^2  \mathcal S_1^{3+L}\phi_{4}\rg]({\bf r}_h^\pm(t))\e^{-i\alpha_nt}\de t.
\end{split}\end{equation}

Similarly, the expression for the corrected approximate scattered field reads as
\begin{equation}\label{eq:approx_scat_fld_3}\begin{split}
u^s_A(\ner)  =&~ (\mathcal D_1^1\phi^c_{A,1})(\ner)-\eta(\mathcal S_{1}^1\phi^c_{A,2})(\ner)+ (\mathcal D_1^2-\gamma\mathcal D_1^3)[w_A\phi^c_{A,3}](\ner)-\\
&~(\mathcal S_1^2-\gamma\mathcal S_1^3)[w_A\phi^c_{A,4}](\ner)+\\
&~\frac{1}{2i}\sum_{n\in\mathcal C_\delta\setminus\mathcal W}\frac{\e^{i\beta_n h}}{\beta_n}\lf\{u_n^-(\ner)\mathsf L_{n}^+\lf[{\bf W}_{\!\!A}\bol\phi_A^c\rg]-u_n^+(\ner) \mathsf L_{n}^-\lf[{\bf W}_{\!\!A}\bol\phi_A^c\rg]\rg\}+\\
&~\frac{1}{2i}\sum_{n\in\mathcal W}\p_{\beta_n}\lf\{u_n^-(\ner)\mathsf L_{n}^+\lf[{\bf W}_{\!\!A}\bol\phi_A^c\rg]-u_n^+(\ner) \mathsf L_{n}^-\lf[{\bf W}_{\!\!A}\bol\phi_A^c\rg]\rg\},\quad \ner\in \Omega'_A.
\end{split}\end{equation}

\subsection{Fredholm property} 

Assuming that $\Gamma_1$ and $\Gamma_2$ are sufficiently smooth, say, with twice continuously differentiable paremetrizations ${\bf r}_1$ and ${\bf r}_2$, respectively,  it is easy to show that the corrected windowed BIE~\eqref{eq:corrected_lin_sys} is Fredholm of the second kind. 

For the sake of presentation simplicity we prove Fredholmness of the corrected windowed BIE~\eqref{eq:corrected_lin_sys}  in the product space $X:=[L^2(0,2\pi)]^2\times [L^2(-A,A)]^2$ for which we first write it as
\begin{equation}\label{eq:win_crr_BIE}
(\mathsf{Id}_X{\bf E}+\bol{\mathsf M}_c\circ\mathsf{Id}_X{\bf W}_A)\bol\phi_A^c= \bol \phi^\inc
\end{equation} where $\bol \phi^\inc\in X$ and the solution $\bol\phi_A^c$ is sought in that same space. Here, $\mathsf{Id}_X$ denotes the identity mapping of $X$ and, slightly abusing the notation, $\bol{\mathsf M}_c$ is considered as an operator acting on $X$, i.e., all the integrals over $\R$ in the definition of $\bol{\mathsf M}_c$  are truncated to the finite interval $[-A,A]$.

Using then the fact that the sub-block operators $\mathsf M_{p,q}$, $p,q=1,\ldots,4,$ defined in~\eqref{eq:blocks1} and~\eqref{eq:blocks2}, are of the Hilbert-Schmidt type (because the associated kernels belong to $L^2([0,2\pi]\times [0,2\pi])$ for $p=q=1,2$, $L^2([0,2\pi]\times [-A,A])$ for $p=1,2$, $q=3,4$, $L^2([-A,A]\times [-A,A])$ for $p=q=3,4$, and $L^2([-A,A]\times[0,2\pi])$ for  $p=3,4$, $q=1,2$) it follows from classical arguments~\cite{atkinson1996numerical} that $\bol{\mathsf M}:X\to X$ is compact. 
On the other hand, since the functionals $\mathsf L_n^\pm,\p_{\beta_n}\mathsf L_n^\pm:X\to \C$ are bounded (because all the integrands involved in their defintion~\eqref{eq:L_functionals_v2} are $L^2$-integrable) and $\Phi_n^\pm,\Phi_n,\p_{\beta_n}\Phi_n^\pm\in X$, we have that the finite-rank operators  $\bol \Phi_n^\pm \mathsf L_n^\pm,\bol\Phi_n\p_{\beta_n}\mathsf L_n^\pm,\p_{\beta_n}\bol\Phi_n\mathsf L_n^\pm:X\to X$ are also compact, and so it is the finite linear combination of them that appears in the definition of $\bol{\mathsf M}_c$ in~\eqref{eq:corrected_mat}. This shows that $\bol{\mathsf M}_c:X\to X$ is compact. 

Therefore, being $\bol{\mathsf M}_c\circ\mathsf{Id}_X{\bf W}_A$ the composition of $\bol{\mathsf M}_c$, which is compact, and $\mathsf{Id}_X{\bf W}_A:X\to X$, which is bounded, we conclude that $\bol{\mathsf M}_c\circ\mathsf{Id}_X{\bf W}_A:X\to X$ is itself compact.  

The Fredholm property of~\eqref{eq:win_crr_BIE} hence follows directly from the invertibility of the operator $\mathsf{Id}_X{\bf E}:X\to X$. 

Having established the Fredholm property of the system~\eqref{eq:corrected_lin_sys}, we can conclude from the Fredholm alternative that existence of solutions in the function space $X$ is implied by uniqueness. We found, however, the uniqueness property difficult to prove since standard arguments based on the unique solvability of associated PDEs (e.g.~\cite{bonnet1994guided}) does not directly apply in this case due to the presence of the windowed integral kernels. Nevertheless, extensive numerical experimentation supports the conjecture that the corrected windowed BIE system~\eqref{eq:corrected_lin_sys} does not suffer from uniqueness issues, which typically manifests at the discrete level as severely ill-conditioned linear systems at certain countable frequencies. 

Finally, we briefly mention that a similar analysis can be carried out in higher-order Sobolev spaces by relying on the well-established mapping properties~\cite{Mclean2000Strongly} of the integral operators~\eqref{eq:operators}.

\begin{figure}[h!]
\centering	
\includegraphics[scale=0.48]{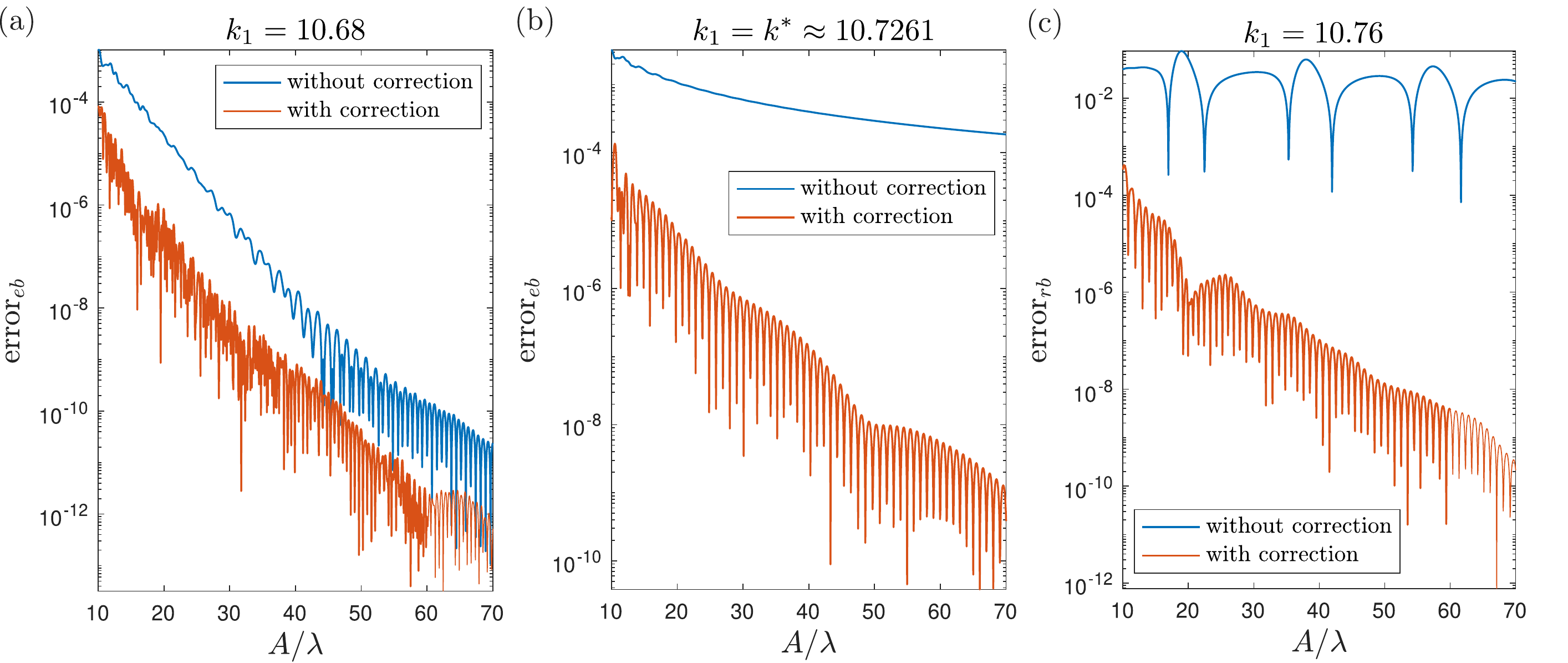}
\caption{Errors~\eqref{eq:error_energy} in the numerical solution of the test problem of Section~\ref{sec:illus_num_exam} obtained using the corrected windowed integral equation~\eqref{eq:corrected_lin_sys}. Three different exterior wavenumbers are considered corresponding to (a)~$k_1=10.68$, (c)~$k_1=10.76$, and (b)~$k_1=k^*$, which corresponds to a RW-anomaly frequency. The fixed parameter value $\delta=3k_1/4$, which yields a four-element set $\mathcal C_\delta$ of correcting terms, is used in all these examples.}\label{fig:example_correction}
\end{figure}
\section{Numerical examples}\label{sec:num_exam}
This section presents a variety of numerical examples that demonstrate the accuracy and robustness of the proposed WGF methodology.
\subsection{Validation examples} 
We start off by applying the proposed windowed BIE approach to the kite-shaped array test problem of Section~\ref{sec:illus_num_exam}, where the naive windowed BIE formulation failed to produce accurate solutions at and around RW-anomaly configurations.

Figure~\ref{fig:example_correction} displays the energy balance errors~\eqref{eq:error_energy} for the problem of scattering by the $2$-periodic array of penetrable kite-shaped obstacle~\eqref{eq:kite_param} for  $k_1\in\{10.68,k^*,10.76\}$ and  $A\in[10\lambda,70\lambda]$ produced by the naive (blue curves) and corrected (red curves) BIE formulations. The same high-order Nystr\"om discretization scheme  was employed to numerically solve both BIEs. The additional parameter $\delta>0$ that enters the corrected BIE~\eqref{eq:corrected_lin_sys} through the  set $\mathcal C_\delta$ in~\eqref{eq:reg_modes}, which selects the modes to be used in the correcting terms in~\eqref{eq:corrected_mat} and~\eqref{eq:approx_scat_fld_3}, is chosen as $\delta = 3k_1/4$ in these examples.  As can be observed in Figure~\ref{fig:example_correction}, the upper envelopes to the red error curves corresponding to  the corrected windowed BIE formulation, exhibit super-algebraic convergence as the window size $A$ increases, for all three wavenumbers considered including the challenging RW-anomaly configuration at $k_1=k^*$.

\begin{figure}[h!]
\centering	
 \includegraphics[scale=0.48]{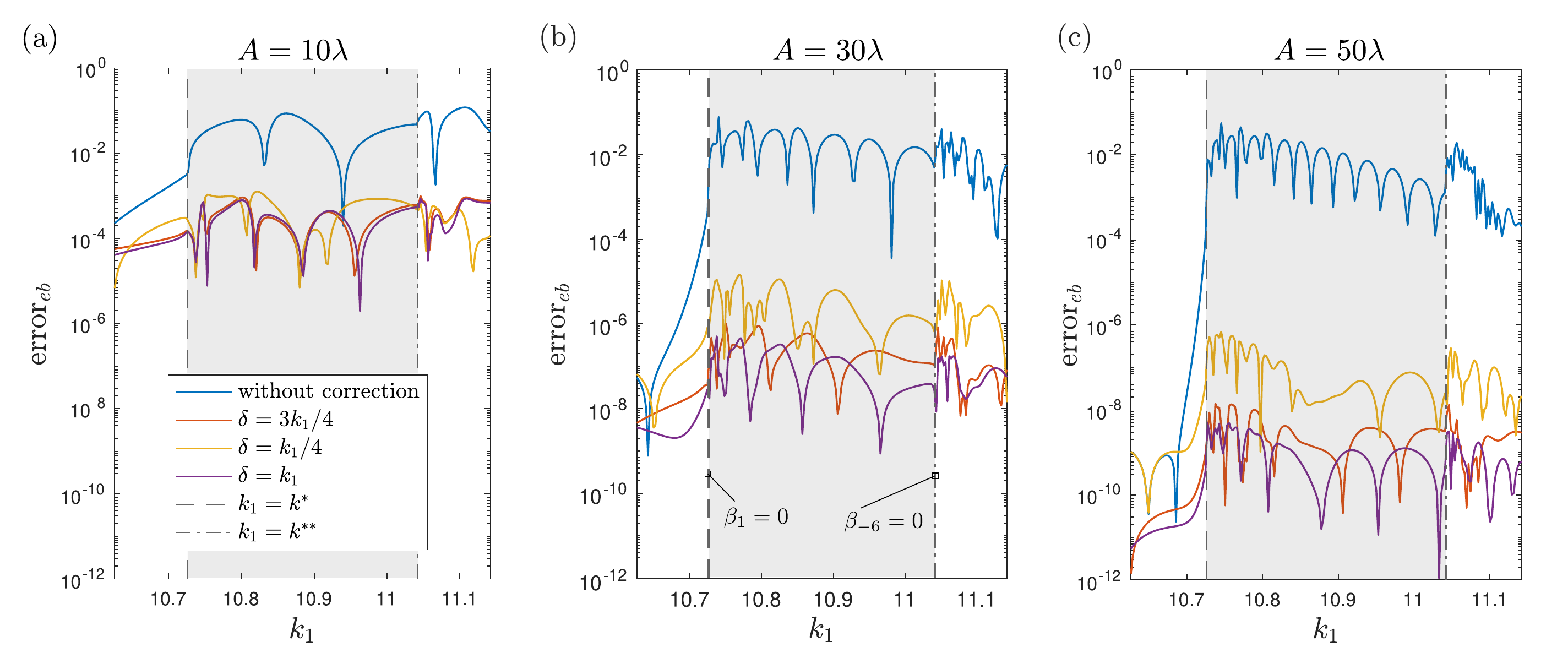}
\caption{Energy balance error~\eqref{eq:error_energy} sweeps for $k_1\in[k^*-0.1,k^{**}+0.1]$, where $k^*$ and $k^{**}$ are two consecutive RW frequencies, in the  solution of the test problem of Section~\ref{sec:illus_num_exam} produced by the corrected windowed BIE~\eqref{eq:corrected_lin_sys} using the parameter values $\delta\in\{k_1/2,3k_1/4,k_1\}$ and (a)~$A=10\lambda$, (b)~$A=30\lambda$, and (c)~$A=50\lambda$.  }\label{fig:sweep_corrected}
\end{figure}

Next, Figure~\ref{fig:sweep_corrected} displays wavenumber sweeps of the energy balance error~\eqref{eq:error_energy} obtained using the naive and the corrected BIE formulations for three window sizes $A\in\{10\lambda,30\lambda,50\lambda\}$ and  $\delta =\{3k_1/4,k_1/4,k_1\}$. The $k_1$-wavenumber range $[k^*-0.1,k^{**}+0.1]$ considered in these examples includes two RW-anomaly frequencies at $k^*\approx10.7261 $ and $k^{**} \approx  11.0418$ where $\beta_1=0$ and $\beta_{-6}=0$, respectively. Unlike the results produced by the naive windowed BIE (blue curves) the corrected approach does not break down at and around RW-anomaly frequencies. Indeed, despite the proximity to the RW frequencies, no extreme accuracy variations are observed as $k_1$ changes while maintaining the main parameters $A$ and $\delta$ fixed. These results demonstrate the robustness of the proposed methodology. Moreover, these results show that the parameter value $\delta=3k_1/4$ (red curves) is good enough to achieve highly accurate solutions throughout the spectrum as no significant improvement is achieved using $\delta=k_1$ (purple curves).

\begin{figure}[h!]
    \centering
    \includegraphics[scale=1]{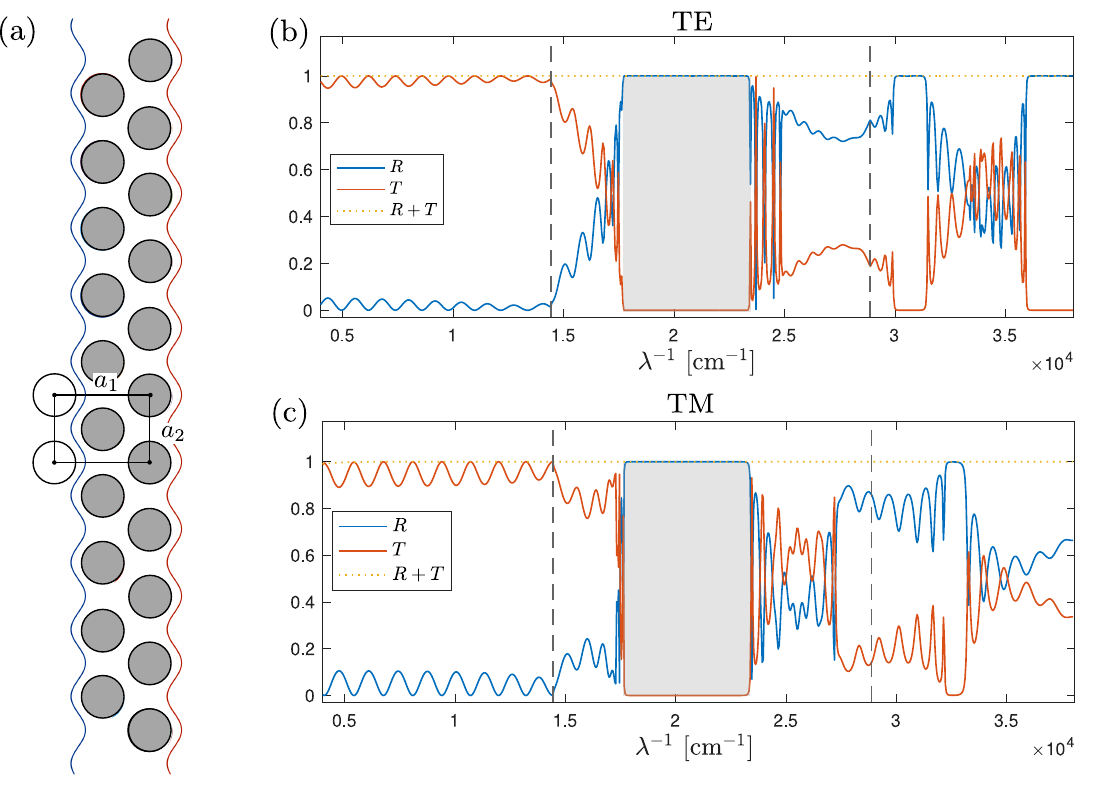}
    \caption{Reflectance and transmittance spectra of a finite-thickness photonic crystal slab in TE and TM polarizations at normal planewave incidence. (a) Depiction of the lattice geometry and the curves involved the numerical solution of the problem by the proposed windowed Green function method. Computed reflectance ($R$) and transmittance ($T$) for various frequencies $\lambda^{-1} = k_1/(2\pi)$ in TE (a) and TM (b) polarization. The first stop band, from $17783\ {\rm cm}^{-1}$ to $23152\ {\rm cm}^{-1}$, is marked in grey, which is the same in both polarizations. The location of RW-anomaly frequencies is marked by the vertical dashed lines.}
    \label{fig:reflectivity}
\end{figure}

\subsection{Photonic crystal slab}\label{sec:PC_example}
In our next and final example we apply the proposed BIE method to the solution of a problem of scattering by a finite-thickness photonic crystal slab. As shown in Figure~\ref{fig:reflectivity}(a) and following the experimental setup of~\cite{huisman2012observation}, we examine a 2D photonic crystal with a centered rectangular lattice of width $a_1=693$~nm and height $a_2=488$~nm.  The refractive index  inside  the crystal--- which is assumed to occupy the exterior domain $\Omega'$---is taken equal to $n=k_1/k_2=2.6$. The boundaries of the~21 pores encompassed by our computational domain (which make up a non-connected curve $\Gamma_1$) are circles of radius $r=155$~nm centered at 
$${\bf a}_l = \frac{(-1)^{l-1}a_1}{4}{\bf e}_1 +\frac{(11-l)a_2}{2}{\bf e}_2, \quad l=1,\ldots 21.$$
  Non-straight unit-cell boundaries $\Gamma_2$ and $\Gamma_3$ parametrized by properly scaled sine functions are used in this example. Note that non-straight curves $\Gamma_2$ and $\Gamma_3$ are necessary in this case to avoid them to intercept the pores ($\Gamma_1$). All the curves involved in the computations are displayed in Figure~\ref{fig:reflectivity}(a) together the lattice geometry. Both TE and TM polarization cases are considered under normal planewave incidence ($\theta^{\inc} = 0$) and the spectrally accurate MK Nystr\"om method is employed in the numerical solution of the corrected windowed BIE~\eqref{eq:win_crr_BIE}.

The computed reflectance ($R$) and transmittance ($T$), which are given by 
\begin{equation}
    R :=\sum_{n\in\mathcal U}\frac{\beta_n}{\beta}\left|B_n^+\right|^2\andtext T :=1+2\real(B_0^-)+\sum_{n\in\mathcal U}\frac{\beta_n}{\beta}\left|B_n^-\right|^2,
\end{equation} are displayed in Figures~\ref{fig:reflectivity}(b)-(c) for TE and TM polarizations, respectively, as functions of the frequency $\lambda^{-1} = k_1/(2\pi)$ in the range from $4000~{\rm cm}^{-1}$ to $38000~{\rm cm}^{-1}$.  Both $R$ and $T$ are here computed using~\eqref{eq:approx_coeff} to approximate the Rayleigh coefficients $B_n^\pm$ and~\eqref{eq:approx_scat_fld_3} to evaluate the scattered field $u^s_A$ on the horizontal lines $y=\pm (5 a_2+2r)$ where coefficients are computed.  The quantity $R+T$,  which is also displayed in those figures, deviates less than $0.01\%$ from its theoretical value of one (see Appendix~\ref{rem:energy}) in all the frequencies considered in this example where we used the parameter values $c=0.5$, $A=20\lambda$, $\delta=3k_1/4$ and $h=cA$ as well as sufficiently refined discretizations of the curves involved. The resulting linear systems, whose sizes remain almost constant around $4920\times 4920$,  were solved by means of GMRES with a tolerance of $10^{-6}$. The observed numbers of GMRES iterations needed to achieved the prescribed tolerance grew with the frequency from 33 (resp. 44) iterations at 4000 ${\rm cm}^{-1}$ to 348 (resp. 470) iterations at 38000 ${\rm cm}^{-1}$ in TE (resp. TM) polarization. The preconditioned system $(\mathsf{Id}_X+\mathsf{Id}_X{\bf E}^{-1}\circ\bol{\mathsf M}_c\circ\mathsf{Id}_X{\bf W}_A)\bol\phi_A^c= \mathsf{Id}_X{\bf E}^{-1}\bol \phi^\inc$  was used in the latter case, as it yields smaller number of iterations.

As expected, band structures form in the reflectance and transmittance spectra displayed in Figures~\ref{fig:reflectivity}(b)-(c). The lowest frequency band structure occurs at roughly the same frequency range in both polarizations between $17783\ {\rm cm}^{-1}$ and $23152\ {\rm cm}^{-1}$ at which $R\approx 1$ in both cases. These results differ slightly from~\cite{huisman2012observation} that places the first stop band for TE-polarized incidence between $4700 n=12220\ {\rm cm}^{-1}$ and $7300 n=18980\ {\rm cm}^{-1}$. 

\begin{figure}[h!]
    \centering
    \includegraphics[scale=1]{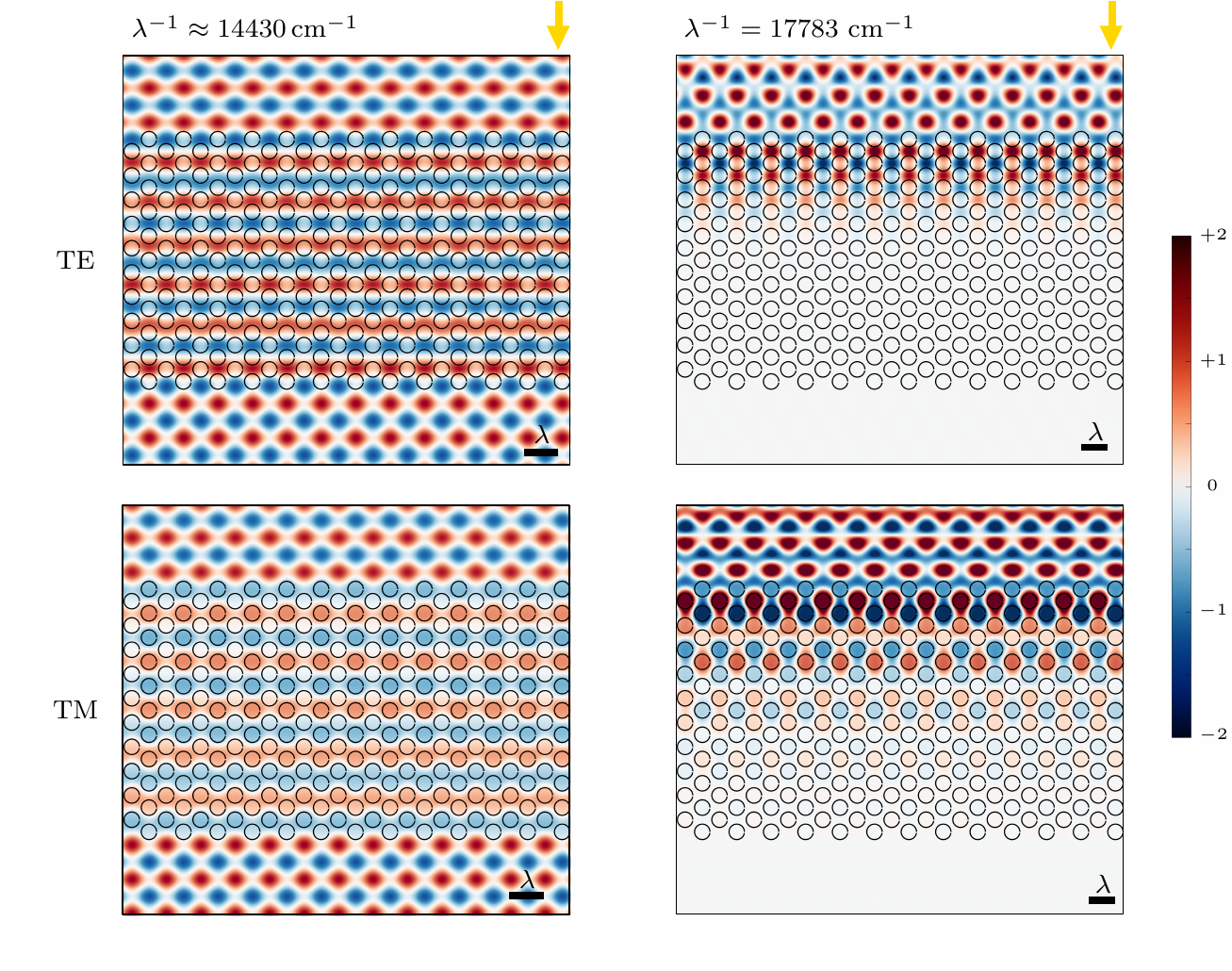}
    \caption{Solution of the problem of scattering of planewave at normal incidence by the finite-thickness photonic crystal of Figure~\ref{fig:reflectivity}(a). Top row: real part of the $z$-component total electric field at the lowest RW-anomaly frequency (left) and at the lowest stop-band frequency (right). Bottom row: real part of the $z$-component of the total magnetic field at the lowest RW-anomaly frequency (left) and at the lowest stop-band frequency (right).}
    \label{fig:PC_plot}
\end{figure}

Finally, Figure~\ref{fig:PC_plot} shows the real part of the total field solution of the problem of scattering by the photonic crystal slab in TE (top row) and TM (bottom row) polarizations for two different frequencies. The left column plots correspond to the RW-anomaly frequency that is marked by the left dashed vertical line in Figures~\ref{fig:reflectivity}(b)-(c). The reflectance $R$ equals $1.5\times 10^{-2}$ and $7.2\times 10^{-4}$ in TE and TM polarization, respectively, at this frequency. The right column plots, on the other hand, correspond to the frequency at the beginning of the stop-band where $R\approx 1$ in both polarizations.

\section{Conclusions and future work}
\label{sec:conclusions}
We presented a novel BIE method for the numerical solution of problems of planewave scattering by periodic line arrays of penetrable obstacles in two dimensions. Our windowed BIE~\eqref{eq:win_crr_BIE}, which is provable Fredholm of the second kind, involves the compact operator $\bol{\mathsf M}_c$~\eqref{eq:corrected_mat} that is expressed in terms of free-space Green function kernels. As such,~\eqref{eq:win_crr_BIE} can be directly discretized and solved by means of any of the various 2D Helmholtz BIE solvers available~(cf.~\cite{hao2014high,klockner2013quadrature,faria2021general,perez2018plane}). We demonstrated through numerical experiments that the combination of our proposed super-algebraically convergent WGF method with the spectrally accurate MK Nystr\"om method, yields a high-order frequency-robust BIE solver that does not break down at and around the challenging RW-anomaly configurations.

This work opens up multiple possible directions for future work. Most of what we presented here applies to Helmholtz scattering problems by line arrays of penetrable obstacles in three dimensions. The extension of this approach to Helmholtz scattering problems by two-dimensional arrays of three-dimensional obstacles is currently being investigated. The main challenge there are the nearly-singular integrals that arise when enforcing the quasi-periodicity condition on the four walls of the unit-cell. This issue, however, can be completely avoided by enforcing such a condition on parts of the boundary of a super-cell. Finally, we mention the natural extension of these results to the full Maxwell's equations in 3D periodic media.

\appendix

\section{Energy conservation principle} \label{rem:energy}
 Energy conservation will be used to assess the accuracy of the proposed boundary integral equation solver. Such principle follows from a direct application of Green second identity over  $R_h = \{(x,y)\in R: |y|<h\}\subset h$ for $h>\max\{h^+,-h^-\}$ is large enough so that $R_h$ contains the obstacle $\Omega$. Applying Green second identity over $R_h\setminus\Omega$ we  obtain 
\begin{equation}\begin{split}
0=&\int_{R_h\setminus\Omega} \{u\Delta\bar u-\bar u\Delta u\}\de\ner=\int_{\p\{R_h\setminus\Omega\}} \left\{u\overline{\p_n v}-\bar u\p_{n}u\right\}\de s\\
=&-2i\operatorname{Im}\left(\int_{\Gamma_1}u\overline{\p_n u}\de s\right) -2i\operatorname{Im}\left(\int_{\Gamma^h_2}u\overline{\p_n u}\de s\right)+2i\operatorname{Im}\left(\int_{\Gamma_3^h}u\overline{\p_{n}u}\de s\right)\\
&+2i\operatorname{Im}\left(\int_{-\frac{L}2}^{\frac{L}2}u(x,h)\overline{\p_yu(x,h)}\de x\right)-2i\operatorname{Im}\left(\int_{-\frac{L}2}^{\frac{L}2}u(x,-h)\overline{\p_y u(x,-h)}\de x\right)
\end{split}\label{eq:GreenExt}\end{equation}
where $\Gamma^h_j=\{(x,y)\in\Gamma_j: |y|<h\}$, $j=2,3$ (see Figure~\ref{fig:green}). Note that, without loss of generality, we have assumed that $\mathsf x_2(t) = -\frac{L}{2}$ for $t\in\R$ such that $\mathsf y_2(t)=h$ and $\mathsf y_2(t)=-h$. 

The integrals over the $\Gamma_2^h$ and $\Gamma_3^h$ in~\eqref{eq:GreenExt} cancel each other by virtue of the fact that
$$
\int_{\Gamma_3^h} u\overline{\p_n u}\de s =\int_{\Gamma_2^h} (\e^{i\alpha L}u)(\overline{\e^{i\alpha L}\p_n u})\de s = \int_{\Gamma_2^h} u\overline{\p_n u}\de s.
$$
Similarly, applying Green second identity inside $\Omega$ and assuming that  $k_2>0$ (i.e., that the medium occupying $\Omega$ is non dissipative), we obtain 
$$
0=\int_{\Omega} \{u\Delta\bar u-\bar u\Delta u\}\de\ner=\int_{\Gamma_1} \left\{u\overline{\p_n u}-\bar u\p_n u\right\}\de s=2i\operatorname{Im}\left(\int_{\Gamma_1} u\overline{\p_n u}\de s\right).
$$
Therefore, we conclude from these identities that
$$
\operatorname{Im}\left(\int_{-\frac{L}2}^{\frac{L}2}u(x,h)\overline{\p_y u(x,h)}\de x\right)=\operatorname{Im}\left(\int_{-\frac{L}2}^{\frac{L}2}u(x,-h)\overline{\p_y u(x,-h)}\de x\right).
$$

These integrals can be expressed in terms of the Rayleigh coefficients. In order to do so, first we note that
\begin{align*}
u(x,\pm h) =&\e^{i\alpha x\mp i\beta h}+\sum_{n\in\mathcal U\cup\mathcal W}B^\pm_n\e^{i(\alpha_n x+\beta_nh)}+\sum_{n\in\mathcal V}B^\pm_n\e^{i\alpha_n x-|\beta_n|h}\\
\overline{\p_y u(x,\pm h)} =&i\beta\e^{-i\alpha x\pm i\beta h}\mp\sum_{n\in\mathcal U\cup\mathcal W} i\beta_n \overline{B^\pm_n}\e^{-i(\alpha_n x+\beta_nh)}\mp\sum_{n\in\mathcal V}|\beta_n| \overline{B^\pm_n}\e^{-i\alpha_n x-|\beta_n|h}.
\end{align*}
Therefore, using the fact that 
\begin{equation}\label{eq:ray_coeff}
\frac{1}{L}\int_{-\frac{L}2}^{\frac{L}2}u(x,\pm h)\e^{-i\alpha_nx}\de x=\begin{cases}
\delta_{0,n}\e^{\mp i\beta h}+B^\pm_{n}\e^{i\beta_nh}&n\in\mathcal U\cup\mathcal W\\
B^\pm_{n}\e^{- |\beta_n|h}&n\in\mathcal V,
\end{cases}
\end{equation}
we arrive at 
\begin{equation}\begin{split}
\frac{1}{L}\int_{-\frac{L}2}^{\frac{L}2}u(x, h)\overline{\p_y u(x, h)}\de x =&i\beta-2\beta\imag (B_0^+\e^{2i\beta h})\\
&-\sum_{n\in\mathcal U\cup\mathcal W}(i\beta_n)|B_n^{+}|^2-\sum_{n\in\mathcal V}|\beta_n||B_n^{+}|^2\e^{-2|\beta_n| h}
\end{split}\label{eq:e1}\end{equation}
\begin{equation}\begin{split}
\frac{1}{L}\int_{-\frac{L}2}^{\frac{L}2}u(x,-h)\overline{\p_yu(x, -h)}\de x =&i\beta+2i\beta \real(B_0^-)\\
+&\sum_{n\in\mathcal U\cup\mathcal W}(i\beta_n)|B_n^{-}|^2+\sum_{n\in\mathcal V}|\beta_n||B_n^{-}|^2\e^{-2|\beta_n| h}.
\end{split}\label{eq:e2}\end{equation}
Finally, taking the imaginary part in both integrals~\eqref{eq:e1} and~\eqref{eq:e2}, and equating them, we obtain the following relation between the Rayleigh coefficients
\begin{equation}\label{eq:conservation_energy}
 0=2\real(B_0^-)+\sum_{n\in\mathcal U}\frac{\beta_n}{\beta}\lf\{|B_n^-|^2+|B_n^+|^2\rg\}
\end{equation}
which expresses the energy conservation principle for this system under consideration.

\section{Super-algebraic decay of windowed oscillatory integrals}\label{rem:sub_al_decay}
 The main argument to establish the  super-algebraic convergence as $A\to\infty$ of the terms in~\eqref{eq:1st_part_v1} corresponding to the propagative modes $\beta_n\in\R_{>0}$ is essentially the repeated use of the integration by parts procedure. In order to illustrate this argument, let us consider the single-layer operator $e_A: ={\mathsf{V}}_1^{1,2}[\chi_{A}^+\e^{i\beta_n\, \cdot}]$ which contributes to the term ${\mathsf M}_{1,4}[\chi_A^+\e^{i\beta_n\, \cdot}]$ where $\mathsf M_{1,4}$ is defined in~\eqref{eq:blocks1}. In detail, we examine the oscillatory integral 
$$
e_A(t) = \frac{i}{4}\int_{cA}^\infty H_0^{(1)}\lf(k_1\sqrt{\lf(\mathsf x_1(t)+\tfrac{L}2\rg)^2+(\mathsf y_1(t)-\tau)^2}\rg)w^c_A(\tau)\e^{i\beta_n\tau}\de \tau,\quad  t\in[0,2\pi).
$$
In view of the addition theorem~\cite{Abramowitz1966Handbook}, i.e., 
\begin{equation}\label{eq:add_them}
H_0^{(1)}\lf(k_1\sqrt{\lf(\mathsf x_1(t)+\tfrac{L}2\rg)^2+(\mathsf y_1(t)-\tau)^2}\rg)= \sum_{\ell=-\infty}^{\infty} H_{\ell}^{(1)}(k_1 |\tau|) J_{\ell}\lf(k_1\varrho(t)\rg) \e^{i \ell (\frac{\pi}2-\vartheta(t))}
\end{equation}
where $\varrho= \sqrt{\lf(\mathsf x_1+\tfrac{L}2\rg)^2+\mathsf y_1^2}$,  $\vartheta=\arctan\big(\frac{\mathsf y_1}{\mathsf x_1+\frac{L}2}\big)$ and $\varrho(t)<cA\leq\left|\tau\right|$, it suffices to estimate the convergence of the integrals
\begin{equation}\label{eq:mode_window_integral}
E^{(n,\ell)}_{A}:= \int_{cA}^\infty H_\ell^{(1)}(k_1\tau)w^c_A(\tau)\e^{i\beta_n \tau} \de \tau, \quad \ell\geq 0,
\end{equation}
as $A\to\infty$.
Performing the change variable $\tau= As$ and letting $\xi(s)=w^c_A(s A)$ where $w^c_A(As)=1-\chi(s,c,1)$ with $\chi$ defined in~\eqref{eq:window_function} (note that it does not depend on $A$), $h_{A,\ell}(s) = \e^{-ik_1As}H_\ell^{(1)}\big(k_1As\big)$ and $\kappa_n =\beta_n+k_1\neq 0$, we arrive at
$$
E^{(n,\ell)}_{A}= A\int_{c}^\infty\xi(s)h_{A,\ell}(s )\e^{i\kappa_n A s}\de s.
$$
Integrating by parts $m>0$ times, the integral above can be recast as
\begin{equation*}\begin{split}
E^{(n,\ell)}_{A}=& \frac{1}{(i\kappa_n)^{m} A^{m-1}}\int_{c}^\infty \e^{i\kappa_n A s}\lf(\frac{\de }{\de s}\rg)^{m}\lf[\xi(s)h_{A,\ell}(s )\rg]\de s\\
=& \frac{1}{(i\kappa_n)^{m} A^{m-1}}\sum_{p=0}^{m}\binom{m}{p}\int_{c}^\infty \e^{i\kappa_n A s}\xi^{(m-p)}(s)h_{A,\ell}^{(p)}(s)\de s
\end{split}\end{equation*}
where we have used  Leibniz’s rule and the fact that $\xi$ together with its derivatives of any order vanish at $s=c$.
 We then conclude that
\begin{equation*}\begin{split}
\lf|E^{(n,\ell)}_{A}\rg|\leq& \frac{\|\xi\|_{C^{m}(\R)}}{|\kappa_n|^{m} A^{m-1}}\lf\{\lf\|h_{A,\ell}\rg\|_{L^1(c,1)}+\sum_{p=1}^{m}\binom{m}{p}\lf\|h_{A,\ell}^{(p)}\rg\|_{L^1(c,\infty)}\rg\}.
\end{split}\end{equation*}
To estimate the $L^1$-norm of $h^{(p)}_{A,\ell}$, $0\leq p\leq m$, we employ a  slight refinement of  Lemma~1 in~\cite{Demanet:2010cc} which  yields 
$$|h_{A,\ell}^{(p)}(s)| = \left|\left(\frac{\mathrm{d}}{\mathrm{d} s}\right)^{p}\left[\mathrm{e}^{-\mathrm{i} k_1As } H_{\ell}^{(1)}(k_1A s)\right]\right| \leq
\frac{1}{\sqrt{8k_1A}}\frac{2^{\ell}P_{p}(\ell)}{|\Gamma(\ell-\frac12)|}s^{-(p+\frac12)}
$$
for  $s \geq c$ and $p,\ell\geq 0$, where $P_p$ are positive-coefficient  polynomials of degree $p$. It hence follows from these bounds that
$$
\lf\|h_{A,\ell}\rg\|_{L^1(c,1)} \leq \frac{1}{\sqrt{8k_1A}}\frac{2^{\ell}P_0(\ell)}{|\Gamma(\ell-\frac12)|}\int_{c}^1 s^{-\frac12}\de s\leq \frac{1}{\sqrt{8k_1A}}\frac{2^{\ell}}{|\Gamma(\ell-\frac12)|}\lf(P_0(\ell)\frac{1- c^{\frac12}}{2}\rg)
$$
and 
\begin{equation*}\begin{split}
 \lf\|h^{(p)}_{A,\ell}\rg\|_{L^1(c,\infty)}\leq&\frac{1}{\sqrt{8k_1A}}\frac{2^{\ell}P_{p}(\ell)}{|\Gamma(\ell-\frac12)|}\int_{c}^\infty s^{-(p+\frac12)} \de s=\\
 &\frac{1}{\sqrt{8k_1A}}\frac{2^{\ell}}{|\Gamma(\ell-\frac12)|}\lf(\frac{P_{p}(\ell)c^{-(p-\frac12)}}{p-\frac12}\rg)
\end{split}\end{equation*}
for $p\geq 1$, and, consequently,
\begin{equation}\label{eq:E_n_bound}
\lf|E^{(n,\ell)}_{A}\rg|\leq \frac{\|\xi\|_{C^{m}(\R)}}{|\kappa_n|^{m} A^{m-1}\sqrt{8k_1A }}\frac{2^{\ell}}{|\Gamma(\ell-\frac12)|}\widetilde P_m(\ell),
\end{equation}
where $\widetilde P_m$ is the $m$-degree polynomial given by
$$\widetilde P_m(\ell) = \lf\{P_0(\ell)\frac{1-c^{\frac12}}{2}+\sum_{p=1}^{m}\binom{m}{p}P_p(\ell)\frac{c^{-(p-\frac12)}}{p-\frac12}\rg\}.$$

With the suitable upper  bounds~\eqref{eq:E_n_bound} for $E_n^{(n,\ell)}$ we return to the addition theorem~\eqref{eq:add_them} to obtain
\begin{equation}\label{eq:sup_alg_conv}\begin{split}
|e_A(t)| &= \frac{1}{4} \lf|\sum_{\ell=-\infty}^{\infty}J_{\ell}\lf(k_1 \varrho(t)\rg) \e^{i \ell (\frac{\pi}2-\vartheta(t))}E^{(n,\ell)}_{A}\rg| \leq \frac{1}{2} \sum_{\ell=0}^{\infty}\lf|J_{\ell}\lf(k_1 \varrho(t)\rg)\rg|\lf|E^{(n,\ell)}_{A}\rg|\\
&\leq  \frac{\|\xi\|_{C^{m}(\R)}}{|\kappa_n|^{m} A^{m-1}\sqrt{32k_1A }}\sum_{\ell=0}^{\infty}a_\ell(t)\quad\text{for all }\quad m\geq 1,\end{split}
\end{equation}
where coefficients in the series above are given by
$a_\ell(t)=\lf|J_{\ell}\lf(k_1 \varrho(t)\rg)\rg|\frac{2^{\ell}\widetilde P_m(\ell)}{|\Gamma(\ell-\frac12)|}.$

Finally, to prove the super-algebraic decay of the function $e_A$ as $A\to\infty$ is suffices to show that the series in~\eqref{eq:sup_alg_conv} converges for all $t\in[0,2\pi)$. In order to do so  we resort to the ratio test. From the asymptotic form of the Bessel functions $J_\ell(x)$ for a fixed real number $x$ and large integers $\ell$, and Stirling's formula~\cite{Abramowitz1966Handbook},
 we readily get that
$$
a_\ell(t) \sim \frac{\widetilde P_m(\ell)}{{2\pi}}\lf(\frac{\e^2k_1\varrho(t)}{\ell(\ell-\frac12)}\rg)^\ell\lf(\frac{\ell}{\e}\rg)^{\frac12} \quad\text{as}\quad \ell\to\infty.
$$
Therefore,
$$
 \lim_{\ell\to\infty} \frac{a_{\ell+1}(t)}{a_\ell(t)} =  \lim_{\ell\to\infty}\lf\{\frac{\e^2k_1\varrho(t)}{(\ell+1)(\ell+\frac12)}\lf(\frac{\ell(\ell-\frac12)}{(\ell+1)(\ell+\frac12)}\rg)^{\ell} \rg\}=0
$$
and hence the desired result follows.

\section*{Acknowledgements}
The authors thank Lars Corbijn van Willenswaard, MACS/COPS, University of Twente, for valuable discussions and for suggesting the finite-thickness photonic crystal example presented in Section~\ref{sec:PC_example}.

\bibliographystyle{abbrv}
\bibliography{References}

\end{document}